\documentclass[12pt,a4paper,reqno]{amsart}
\usepackage{amsmath}
\usepackage{amssymb,latexsym}
\usepackage[dvips]{graphicx}
\usepackage{color}
\usepackage{cite}
\usepackage{amsthm}

\newtheorem{theorem}{Theorem} [section]
\newtheorem{lemma}{Lemma} [section]
\newtheorem{proposition}{Proposition} [section]
\newtheorem{corollary}{Corollary} [section]

\newtheorem{remark}{Remark}[section]

\let\ssection=\section\renewcommand{\section}{\setcounter{equation}{0}\ssection}
\begin{document}
\date{}

\address{M. Darwich: Universit\'e Fran\c{c}ois rabelais de Tours, Laboratoire de Math\'ematiques
et Physique Th\'eorique, UMR-CNRS 6083, Parc de Grandmont, 37200
Tours, France} \email{Mohamad.Darwich@lmpt.univ-tours.fr}

\title[Blowup]{BLOWUP FOR THE DAMPED $L^{2}$-CRITICAL NONLINEAR SCHR\"{O}DINGER EQUATION}
\author{Darwich Mohamad}

\begin{abstract}
We consider the Cauchy problem for the  $L^{2}$-critical damped nonlinear Schr\"{o}dinger equation. We prove existence and stability of finite time blowup dynamics with the log-log blow-up speed  for $\|\nabla u(t)\|_{L^2}$.
\end{abstract}

\maketitle
%%%%%%%%%%%%%%%%%%%%%%%%%%%%%%%%%%%%%%%%%%%%%%%%%%%%%%%%%%%%%%%%%%%%%%%%%%%%%%%%%%%%%%%%%%%%%%%%%%%%%%%%%%%%%%%%%%%%%%%%%%%%%%%%%%%
\section{Introduction}
In this paper, we study the blowup of  solutions to the Cauchy problem for  the $L^{2}$-critical damped nonlinear Schr\"{o}dinger equations:
\begin{equation}\label{NLSa}
\begin{cases}
%$$ \left\lbrace \begin{array}{l}
iu_{t} + \Delta{u} +|u|^{\frac{4}{d}}u + iau =0,  (t,x) \in [0,\infty[\times \mathbb{R}^{d}, d=1,2,3,4. \\
u(0)= u_{0} \in H^1(\mathbb{R}^{d})
%\end{array} $$
\end{cases}
\end{equation}
with initial data $u(0)= u_{0} \in H^1(\mathbb{R}^{d})$ and where $a > 0$ is the coefficient of friction. Equation (\ref{NLSa}) arises in various areas of nonlinear optics, plasma physics and fluid mechanics.
It is known that the Cauchy problem for (1.1) is locally well-posed in $H^1(\mathbb{R}^{d})$(see Kato\cite{Kato} and also Cazenave\cite{Cazenave}): For any $u_{0} \in H^{1}(\mathbb{R}^{d})$, there exist $T \in (0,\infty]$ and a unique solution $u(t)$ of $(1.1)$ with $u(0)=u_{0}$ such that $u \in C([0,T);H^1(\mathbb{R}^{d}))$. Moreover, T is the maximal existence time of the solution $u(t)$ in the sense that if $T < \infty$ then $\displaystyle{ \lim_{t\rightarrow T}{\|u(t)\|_{H^1(\mathbb{R}^{d})}}}=\infty$.\\
Ohta \cite{ohta} and Tsutsumi \cite{Tsutsumi} studied the supercritical case($|u|^{p}u$ with $p > \frac{4}{d}$) and showed that blow-up in finite time can occur, using the virial method. However this method does not seem to apply in the critical case. Therefore, even if numerical simulations suggest the existence of finite time blowup solutions in this case(see Fibich \cite{Fibich}), there does not exist any mathematical proof of blow-up in the critical case.\\
Let us notice that for $a=0$ (\ref{NLSa}) becomes the $L^2$-critical nonlinear Schr\"{o}dinger equation:\\
\begin{equation}\label{NLS}
\begin{cases}
%$$ \left\lbrace \begin{array}{l}
 iu_{t} + \Delta u + |u|^{\frac{4}{d}}u = 0\\
 u(0)=u_{0} \in H^{1}(\mathbb{R}^{d})
 %\end{array} $$
 \end{cases}
 \end{equation}
 This equation (\ref{NLS}) admits a number of symmetries in the energy space $H^1$: if $u(t, x)$ is a
solution to (\ref{NLS}) then $\forall \lambda_{0} \in \mathbb R$, so is $\lambda_{0}^{\frac{d}{2}}u(\lambda_{0}x,\lambda_{0}^{2}t)$.
Note that the $L^2$-norm is left invariant by the the scaling symmetry and thus $L^2$ is the critical
space associated with this symmetry.\\
The evolution of (\ref{NLS}) admits the following conservation laws in the energy space $H^1$:\\
$L^2$norm : $\left\|u(t, x)\right\|_{L^2} = \left\|u(0, x)\right\|_{L^2} = \left\|u_{0}(x)\right\|_{L^2}.$\\
Energy : $E(u(t, x)) = \frac{1}{2}\|\nabla u\|_{L^2}^{2} - \frac{d}{4 + 2d}\|u\|_{L^{\frac{4}{d}+2}}^{\frac{4}{d}+2} = E(u_{0}).$\\
Kinetic momentum : $P(u(t))=Im(\displaystyle{\int} \nabla u \overline{u}(t,x))= P(u_{0}).$\\
Special solutions play a fundamental role for the description of the dynamics of (\ref{NLS}). They are the solitary waves of the form $u(t, x) =\exp(it)Q(x)$, where $Q$ solves:
 \begin{equation}\label{ellip}
 \Delta Q + Q|Q|^{\frac{4}{d}} = Q.
 \end{equation}
 Equation (\ref{ellip}) is a standard nonlinear elliptic equation,  that possesses  a unique positive solution (see \cite{Ber}, \cite{Lions1},\cite{Kwong1}) .\\
For $u_0 \in  H^1$, a sharp criterion for global existence has been exhibited by Weinstein
\cite{weinst}:\\
a) For $\|u_0\|_{L^2} < \|Q\|_{L^2}$ , the solution of (\ref{NLS}) is global in $H^1$. This follows from the conservation
of the energy and the $L^2$ norm and the sharp Gagliardo-Nirenberg inequality:
\begin{equation}
\forall u \in H^1,  E(u) \geq \frac{1}{2}(\int |\nabla u|^2)\bigg(1 - \big(\frac{\int |u|^2}{\int |Q|^2}\big)^{\frac{2}{d}}\bigg). \nonumber
\end{equation}
b)There exists blow-up solutions emanating from initial data $u_0 \in H^1$ with $\left\|u_0\right\|_{L^2} \leq \left\|Q\right\|_{L^2}$. This follows from the pseudo-conformal symmetry applied to the solitary waves.
 In the series of papers \cite{MerleRaphael1,Merle6}, Merle and Raphael have studied the blowup for the $L^{2}$-critical nonlinear Schr\"odinger equation (\ref{NLS}) and have proven the  existence of the blowup regime corresponding to the log-log law:
\begin{equation}\label{speed}
\displaystyle{\|u(t)\|_{H^1(\mathbb{R}^{d})} \sim \bigg(\frac{\text{log}\left|\text{log}(T-t)\right|}{T-t}\bigg)^{\frac{1}{2}}.}
\end{equation}
This regime has the advantage to be stable with respect to $H^1$-perturbation and with respect some perturbations of the equation.

\begin{remark}\label{remark}
 Based on the works \cite{MerleRaphael1,Merle5} we have the following result:\\
Let $u_{0}$  the initial data $\in H^{1}(\mathbb{R}^{d})$ with small super-critical mass:
\begin{equation}\label{alpha}
\|Q\|_{L^2} < \|u_{0}\|_{L^2} < \|Q\|_{L^2} + \alpha_{0}
\end{equation}\\
with nonpositive Hamiltonian $E(u_{0}) < 0$, then the corresponding solution to (\ref{NLS}) blowup in finite time with the log-log speed.

\end{remark}
In the case of (\ref{NLSa}), there does not exists conserved quantities anymore. However, it is easy to prove that if $u$ is a solution of (\ref{NLSa}) then:
 \begin{equation}\label{mass a}
 \|u(t)\|_{L^2}=\exp (-at)\|u_{0}\|_{L^2}, t \in [0,T),
 \end{equation}
 \begin{equation}\label{derivee de lenergie}
 \frac{d}{dt}E(u(t))=-a(\|\nabla u\|_{L^{2}}^{2} - \|u\|_{L^{\frac{4}{d}+2}}^{\frac{4}{d}+2})
 \end{equation}
 and
 \begin{equation}\label{moment}
 |P(u(t))|=\exp (-2at)|P(u_{0})|, t \in [0,T).
 \end{equation}
In this paper, we will show that:
\begin{enumerate}
\item if $\|u_0\|_{L^2} \leq \|Q\|_{L^2}$, then the solution of (\ref{NLSa}) is global in $H^1$.
\item The existence of finite time blowup solutions. 
\end{enumerate}
More precisely, we have the following theorem:
\begin{theorem}\label{theoremessentiel}
 Let $u_0$ in $H^1(\mathbb{R}^{d})$ with $d=1,2,3,4$:
\begin{enumerate}
 \item if \quad $\|u_0\|_{L^2} \leq \|Q\|_{L^2}$ then the solution of (\ref{NLSa}) is global in $H^1$.
\item There exists $\delta_0 > 0$ such that $\forall a > 0$ and $\forall \delta \in]0, \delta_{0}[$, there exists $u_0 \in H^1$ with $\|u_0\|_{L^2} = \|Q\|_{L^2} + \delta$, such that the solution of (\ref{NLSa}) blows up in finite time in the log-log regime.
\end{enumerate}
\end{theorem}
To show the existence of the explosive solutions, we will put us  in the log-log regime described by Merle and Raphael.\\
The global existence will be proved thanks to a $L^2$-concentration phenomenon (see Proposition \ref{nonradiale} in the next section).\\

$\textbf{Acknowledgments.}$ I  would like  to thank prof Luc Molinet for his encouragement,
advice, help and for the rigorous attention to this paper.\\

\section{$L^2$-concentration }
In this section, we prove assertion (1) of Theorem \ref{theoremessentiel} by extending the proof of the $L^2$-concentration phenomen, proved by Ohta and Todorova \cite{ohta} in the radial case, to the non radial case.\\
 Hmidi and Keraani  showed in \cite{Hmidi} the  $L^2$-concentration for the equation (\ref{NLS}) without the hypothese of radiality, using the following theorem:
 \begin{theorem}\label{limsup}
 Let $(v_{n})_{n}$ be a bounded family of $H^1(\mathbb{R}^{d})$, such that:
 \begin{equation}
 \limsup_{n \rightarrow +\infty}\left\|\nabla v_{n}\right\|_{L^2(\mathbb{R}^{d})} \leq M \quad and \quad \limsup_{n \rightarrow +\infty}\left\|v_{n}\right\|_{L^{\frac{4}{d} + 2}} \geq m.
 \end{equation}
 Then, there exists $(x_{n})_{n} \subset \mathbb{R}^{d}$ such that:
 \begin{equation}
 v_{n}(\cdot + x_{n}) \rightharpoonup V \quad weakly, \nonumber
 \end{equation}
 with $\left\|V\right\|_{L^2(\mathbb{R}^{d})} \geq (\frac{d}{d+4})^{\frac{d}{4}}\frac{m^{\frac{d}{2}+1} + 1}{M^{\frac{d}{2}}}\left\|Q\right\|_{L^2(\mathbb{R}^{d})}$.
 \end{theorem}
 Now we have the following theorem:
 \begin{theorem}\label{nonradiale}
 Assume that $u_{0} \in H^{1}(\mathbb{R}^{d})$ , and suppose that the solution of (\ref{NLSa}) with $u(0)=u_{0}$ blows up in finite time $T \in (0,+\infty)$. Then, for any function $w(t)$ satisfying $w(t)\left\|\nabla u(t)\right\|_{L^2(\mathbb{R}^{d})} \rightarrow \infty$ as $t \rightarrow T$, there exists $x(t) \in \mathbb{R}^{d}$ such that, up to a subsequence,
 \begin{equation}
 \displaystyle{\limsup_{t \rightarrow T}\left\|u(t)\right\|_{L^2(\left|x - x(t)\right| < w(t))} \geq \left\|Q\right\|_{L^2(\mathbb{R}^{d})}.}\nonumber
 \end{equation}
 \end{theorem}
 To show this theorem we shall need the following lemma:
 \begin{lemma}\label{lemma ohta}
  Let $T \in (0,+\infty)$, and assume that a function $F : [0, T )\longmapsto(0,�+\infty)$ is
continuous, and $\lim_{t \rightarrow T} F(t) = +\infty$. Then, there exists a sequence $(t_{k})_{k}$ such that
$t_{k}\rightarrow T $and
\begin{equation}\label{Fsur son int}
\displaystyle{\lim_{t_k \rightarrow T}\frac{\displaystyle{\int}_{0}^{t_{k}}F(\tau)d\tau}{F(t_{k})} = 0.}
\end{equation}
\end{lemma}
For the proof see \cite{ohta}.\\

\textbf{Proof of Theorem \ref{nonradiale}}:\\
By the energy identity $(\ref{derivee de lenergie})$, we have 
\begin{equation}\label{energie}
\displaystyle{E(u(t)) = E(u_{0}) - a\int_{0}^{t}K(u(\tau))d\tau, \quad t \in [0,T[.}
\end{equation}
Where $K(u(t)) = \|\nabla u\|_{L^{2}}^{2} - \|v\|_{L^{\frac{4}{d}+2}}^{\frac{4}{d}+2}$, and by the Gagliardo-Nirenberg inequality and (\ref{mass a}), we have:
\begin{align}
\left|K(u(t))\right| & \leq \left\|\nabla u(t)\right\|_{L^2(\mathbb{R}^{d})}^{2} + \left\|u(t)\right\|_{L^{2 + \frac{4}{d}}}^{2 + \frac{4}{d}}\nonumber\\
&\leq  \left\|\nabla u(t)\right\|_{L^2(\mathbb{R}^{d})}^{2} + C\left\|u(t)\right\|_{L^2(\mathbb{R}^{d})}^{\frac{4}{d}}\left\|\nabla u(t)\right\|_{L^2(\mathbb{R}^{d})}^{2}\nonumber\\
&\leq(1 + C\left\|u_0\right\|_{L^2(\mathbb{R}^{d})}^{\frac{4}{d}})\left\|\nabla u(t)\right\|_{L^2(\mathbb{R}^{d})}^{2}\nonumber
\end{align}
for all $t \in [0,T[$. Moreover, we have $\displaystyle{\lim_{t \rightarrow T}\left\|\nabla u(t)\right\|_{L^2(\mathbb{R}^{d})}} = +\infty$, thus by Lemma \ref{lemma ohta}, there exists a sequence $(t_{k})_{k}$ such that $t_{k} \rightarrow T$ and 
\begin{equation}\label{ksurnabla}
\displaystyle{\lim_{k \rightarrow \infty}\frac{\displaystyle{\int}_{0}^{t_{k}}K(u(\tau))d\tau}{\left\|\nabla u(t_k)\right\|_{L^2(\mathbb{R}^{d})}^{2}} = 0.}
\end{equation}
Let
 $$\rho(t) = \frac{\left\|\nabla Q\right\|_{L^2(\mathbb{R}^{d})}}{\left\|\nabla u(t)\right\|_{L^2(\mathbb{R}^{d})}} \quad \text{and} \quad v(t,x)=\rho^{\frac{d}{2}}u(t,\rho x)$$
and $\rho_{k} = \rho(t_{k}), v_{k} = v(t_{k},.)$. The family $(v_{k})_{k}$ satisfies 

$$\left\|v_{k}\right\|_{L^2(\mathbb{R}^{d})} \leq \left\|u_{0}\right\|_{L^2(\mathbb{R}^{d})}\quad \text{and} \quad \left\|\nabla v_{k}\right\|_{L^2(\mathbb{R}^{d})} = \left\|\nabla Q\right\|_{L^2(\mathbb{R}^{d})}.$$
\bigskip
By (\ref{energie}) and (\ref{ksurnabla}), we have
\begin{equation}\label{Edevk}
\displaystyle{E(v_{k}) = \rho^2_{k}E(u_{0}) -a\rho^2_{k}\int_{0}^{t_{k}}K(u(\tau))d\tau \rightarrow 0,}
\end{equation}
which yields 
\begin{equation}\label{vk ver Q}
\displaystyle{\left\|v_{k}\right\|_{L^{\frac{4}{d} + 2}}^{\frac{4}{d} + 2} \rightarrow \frac{d + 2}{d}\left\|\nabla Q\right\|_{L^2(\mathbb{R}^{d})}^{2}.}
\end{equation}
The family $(v_{k})_{k}$ satisfies the hypotheses of Theorem \ref{limsup} with \\
$$m^{\frac{4}{d} + 2} = \frac{d+2}{d}\left\|\nabla Q\right\|_{L^2(\mathbb{R}^{d})}^{2} \quad \text{and} \quad M = 
\left\|\nabla Q\right\|_{L^2(\mathbb{R}^{d})},$$
\bigskip
thus there exists a family $(x_{k})_{k} \subset \mathbb{R}^{d}$ and a profile $V \in H^{1}(\mathbb{R}^{d})$ with $\left\|V\right\|_{L^2(\mathbb{R}^{d})} \geq \left\|Q\right\|_{L^2(\mathbb{R}^{d})}$, such that,
\begin{equation}\label{convergencefaible}
\displaystyle{\rho^{\frac{d}{2}}_{k}u(t_{k}, \rho_{k}\cdot +  x_{k}) \rightharpoonup V \in H^{1} \quad \text{weakly}.}
\end{equation}
Using (\ref{convergencefaible}), $\forall A \geq 0$
\begin{equation}
 \displaystyle{\liminf_{n\to +\infty}\int_{B(0,A)}\rho_{n}^{d}|u(t_{n},\rho_{n}x+x_{n})|^{2}dx\geq \int_{B(0,A)}|V|^{2}dx,}\nonumber
 \end{equation}
  but $\lim_{n\to +\infty}\frac{w(t_{n})}{\rho_{n}}=+\infty$\,\,thus $\frac{w(t_{n})}{\rho_{n}}> A$, $\rho_{n}A < w(t_{n})$. This gives immediately:
  
  \begin{align}
  \displaystyle{\liminf_{n\to +\infty}\sup_{y\in\mathbb{R}^{d}}\int_{|x-y|\leq w({t_{n})}}|u(t_{n},x)|^{2}dx\geq \int_{|x|\leq A}|V|^{2}dx.}\nonumber
  \end{align}
  This it is true for all $A > 0$ thus :
  
  \begin{equation}
  \displaystyle{\liminf_{t\to T}\sup_{y\in\mathbb{R}^{d}}\int_{|x-y|\leq w(t)}|u(t,x)|^{2}dx\geq \int Q^{2}.}
  \end{equation}
  
But for every  $t \in[0,T[$, $y \mapsto \displaystyle{\int}_{|x-y|\leq w(t)}|u(t,x)|^{2}dx$ is continuous
and goes to $0$ at infinity, thus the sup is reached in a  point $x(t)\in\mathbb{R}^{d}$, $\displaystyle{\sup_{y\in\mathbb{R}^{d}}\int_{|x-y|\leq w(t)}|u(t,x)|^{2}dx= \int_{|x-x(t)|\leq w(t)}|u(t,x)|^{2}dx}$ and Theorem \ref{nonradiale} is proved.\\
Now the part one of Theorem \ref{theoremessentiel} is a consequence of Theorem \ref{nonradiale} and (\ref{mass a}).\\

\section{ Strategy of the proof of Theorem \ref{theoremessentiel} part 2.}

We look for a solution of (\ref{NLSa}) such that for $t$ close enough to blowup time, we shall have the following decomposition:
\begin{equation}\label{decomposition}
u(t,x)=\frac{1}{\lambda^{\frac{d}{2}}(t)}(Q_{b(t)} + \epsilon)(t,\frac{x-x(t)}{\lambda(t)})e^{i\gamma(t)},
\end{equation}

for some geometrical parameters $(b(t),\lambda(t), x(t),\gamma(t)) \in (0,\infty)\times(0,\infty)\times\mathbb{R}^{d}\times\mathbb{R}$, here $\lambda(t)\sim \frac{1}{\|\nabla u(t)\|_{L^2}}$, and the profiles $Q_{b}$ are suitable deformations of $Q$ related to some extra degeneracy
of the problem.\\

Now we take $u_{0}$ in $H^{1}$  such that $u_{0}$  admits the following controls:\\
\begin{enumerate}
\item Control of the scaling parameter:
\begin{equation}\label{condition sur b et lambda}
0 < b(0) \ll 1 \quad \text{and} \quad 0 < \lambda(0) < e^{-{e^{\frac{2\pi}{3b(0)}}}}.
\end{equation}
\item $L^2$ control of the excess of mass:
\begin{equation}\label{small norm}
\|\epsilon(0)\|_{L^{2}} \ll 1.
\end{equation}
\item $H^1$ smallness of $\epsilon(0)$:
\begin{equation}\label{smallness dans H1}
\displaystyle{\int |\nabla \epsilon(0)|^2 + \int |\epsilon(0)|^2e^{-|y|} \leq \Gamma_{b(0)}^{\frac{3}{4}}.}
\end{equation}
\item Control of the energy and momentum:
\begin{equation}\label{control de l'energie}
|E(u_{0})| \leq \frac{1}{\sqrt{\lambda(0)}}
\end{equation}
\begin{equation}\label{control de moment}
|P(u_{0})| \leq \frac{1}{\sqrt{\lambda(0)}}.
\end{equation}
\end{enumerate}
\begin{remark}
To prove that there exists $u_{0}$ in $H^1(\mathbb{R}^{d})$ satisfying (\ref{condition sur b et lambda})-(\ref{control de moment}), we take  $\tilde{u}_{0}$ in $H^1(\mathbb{R}^{d})$ an initial data such that the corresponding solution to (\ref{NLS}) blows up in the log-log
regime as described by Merle and Raphael. Then from \cite{Merle4} there exists a time $t_{0}$ such that
$\tilde{u}(t_{0})$ admits a geometrical decomposition:
\begin{equation}
\displaystyle{\tilde{u}(t_{0},x)=\frac{1}{\lambda(t_{0})}\big(Q_{b(t_{0})} + \epsilon(t_{0})\big)\bigg(\frac{x-x(t_{0})}{\lambda(t_{0})}\bigg)e^{i\gamma(t)}}\nonumber
\end{equation}
such that (\ref{condition sur b et lambda})-(\ref{smallness dans H1}) hold. Moreover by conservation of the Hamiltonian and the Kinetic momentum:
\begin{equation}
\displaystyle{\left|E(\tilde{u}(t_{0}))\right| + \left|P(\tilde{u}(t_{0}))\right| = \left|E(\tilde{u}_{0})\right| + \left|P(\tilde{u}_{0})\right| \leq \frac{1}{\sqrt{\left\|\nabla \tilde{u}(t_0)\right\|_{L^2}}}}\nonumber
\end{equation}
for $t_{0}$ close enough to blowup time, and hence (\ref{control de l'energie}) and (\ref{control de moment}) hold. We take $u_{0} = \tilde{u}(t_{0})$.
\end{remark}
These conditions will be denoted by $C.I$. Now we have the following theorem:

\begin{theorem}\label{theorem 3}
 Let  $u_{0} \in H^1$ satisfying C.I, then for $0 < a < a_0$, $a_0$ small the corresponding solution $u(t)$ of (\ref{NLSa}) blows up in finite time in the log-log regime.
\end{theorem}
The set of initial data satisfying C.I is open in $H^1$, using the continuity  with regard to the initial data and the parameters, we can prove the following corollary(see the proof in section 5):
\begin{corollary}\label{colloraireperturbation}
Let $u_{0} \in H^1$ be an initial data  such that the corresponding solution $u(t)$ of (\ref{NLS}) blows up in the loglog regime. There exist $\beta_{0} > 0$ and $a_{0} > 0 $ such that if $v_{0} = u_{0} + h_{0}$, $\left\|h_{0}\right\|_{H^{1}} \leq \beta_{0}$ and $a \leq a_{0}$, the solution $v(t)$ for (\ref{NLSa}) with the initial data $v_{0}$ blowup in finite time.
\end{corollary}
\begin{remark}
\begin{enumerate}
\item Combining Theorem \ref{theorem 3} and Remark \ref{remark}, we deduce that for $u_{0} ~ \in ~ H^1(\mathbb{R}^{d})$, having negative energy and satisfying (\ref{alpha}), there exists $a_0 > 0$ such that the corresponding solution of (\ref{NLSa}) blows up in finite time providing $a \leq a_0$.
\item It is easy to check that if $u$ is a solution of $\text{NLS}_{a}$ (Equation (\ref{NLSa})),then $\lambda^{\frac{d}{2}}u(\lambda^2t,\lambda x)$ is a solution of $\text{NLS}_{\lambda^2a}$. Therefore Corollary  \ref{colloraireperturbation} ensures the existence for any $ a > 0 $ of explosive solutions emanating from an initial data $u_{0,a} \in H^1$, where $u_{0,a}$ satisfies:\\
\centerline{$\|Q\|_{L^2} < \left\|u_{0,a}\right\|_{L^2} = \left\|u_{0}\right\|_{L^2} < \|Q\|_{L^2} + \alpha_0$.}
\end{enumerate}
\end{remark} 

%\subsection{Strategy of the proof for theorem \ref{theoremessentiel} part 2}

After the decomposition (\ref{decomposition}) of $u$, the
log-log regime corresponds to the following asymptotic controls
\begin{equation}\label{stragerie 1}
b_{s}\sim Ce^{-\frac{c}{b}},\, -\frac{\lambda_{s}}{\lambda} \sim b
\end{equation}
and 
\begin{equation}\label{stragerie 2}
\int|\nabla \epsilon|^2\lesssim e^{-\frac{c}{b}},
\end{equation}
where we have introduced the rescaled time $\frac{ds}{dt} = \frac{1}{\lambda^2}$.

In fact, (\ref{stragerie 2}) is partly a consequence of the
preliminary estimate:
\begin{equation}\label{stragerie 3}
\int|\nabla \epsilon|^2 \lesssim e^{-\frac{c}{b}} + \lambda^2 E(t).
\end{equation}
One then observes that in the log-log regime, the integration of the laws (\ref{stragerie 1}) yields
\begin{equation}\label{stragerie 4}
\lambda \sim e^{-e^{\frac{c}{b}}} \ll e^{-{\frac{c}{b}}},
\, b(t) \rightarrow 0, \, t\rightarrow T.
\end{equation}
Hence, the term involving the conserved Hamiltonian is asymptotically negligible
with respect to the leading order term $e^{-{\frac{c}{b}}}$
which drives the decay $(\ref{stragerie 3})$ of $b$. This
was a central observation made by Planchon and Raphael in \cite{Planchon 1}. In fact, any growth
of the Hamiltonian algebraically below $\frac{1}{\lambda^{2}}$ would be enough.
In this paper, we will prove that in the log-log regime, the growth of the energy is estimated by:
\begin{equation}\label{decroissance}
\displaystyle{E(u(t)) \lesssim \big(\text{log}(\lambda (t))\big)^2.}
\end{equation}
%where $\alpha > 0$.
We deduce from $(\ref{stragerie 3})$ that:
\begin{equation}
\int|\nabla \epsilon|^2  \lesssim e^{-{\frac{c}{b}}}.
\end{equation}
An important feature of this estimate of $H^1$ flavor is that it relies on a flux computation
in $L^2$ . This allows one to
recover the asymptotic laws for the geometrical parameters $(\ref{stragerie 1})$ and to close the
bootstrap estimates of the log-log regime.\\
   This paper is organized as follows. In Sect. 3, 4 and 5, we recall some nonlinear objects
involved in the $H^1$ description of the log-log regime and set up the bootstrap argument, see Proposition \ref{bootstrap lemma}. In Sect. 6, we will control in the bootstrap regime
the growth of the energy and momentum, see Lemma \ref{control energie}. In
Sect. 7 and 8, we close the bootstrap estimates and conclude the proof of Theorem \ref{theoremessentiel} (part 2). Finally, the 
 proof of Corollary \ref{colloraireperturbation}
  is postponed at the end of Section 8.

\section{Choice of the blow up profile}
Let us introduce the rescaled time :
\begin{equation}
\displaystyle{s(t) = \int_{0}^{t}\frac{d\alpha}{\lambda^2(\alpha)}}.\nonumber
\end{equation}
It is elementary to check that, whatever the behavior of $u(t)$ is, one always has:\\
\centerline{$s([0,T[)=\mathbb{{R}^{+}}.$}
Let us set:\\
\centerline{$v(s,y)=e^{-i\gamma(t)}\lambda^{\frac{d}{2}}u(t,\lambda(t)x + x(t))$,}
where $y = \lambda(t)x + x(t)$, note that:
$$v_s = -i\gamma_s v + \frac{d}{2}\frac{\lambda_s}{\lambda}v + e^{-i\gamma(t)}\lambda^{2 + \frac{d}{2}}u_t + \frac{\lambda_s}{\lambda}y\cdot\nabla v + \frac{x_s}{\lambda}\cdot\nabla v.$$
$$\Delta v =  e^{-i\gamma(t)}\lambda^{2 + \frac{d}{2}}\Delta u(t,\lambda(t)x + x(t))~~ \text{and}~~ v|v|^{\frac{4}{d}} = e^{-i\gamma(t)}\lambda^{2 + \frac{d}{2}}u|u|^{\frac{4}{d}}.$$
Now  $u(t,x)$ solves (\ref{NLSa}) on $[0,T[$ iff $v(s,y)$ solves: $\forall s\geq 0$,
\begin{equation}\label{nouvel equation}
iv_{s} + \Delta v - v + v|v|^{\frac{4}{d}} = i\frac{\lambda_{s}}{\lambda}(\frac{d}{2}v + y\cdot\nabla v) + i\frac{x_{s}}{\lambda}\cdot\nabla v + \tilde{\gamma_{s}}v,
\end{equation}
where $\tilde{\gamma_{s}} = -\gamma_{s} -1 - ia\lambda^{2}$, and $a$ is the coefficient of friction. Now $v(s,y)= Q(y) + \epsilon(s,y)$ and we linearize (\ref{nouvel equation}) close to Q. The obtained system has the form:
\begin{equation}\label{equation en epsilon}
i\epsilon_{s} + L\epsilon = i\frac{\lambda_{s}}{\lambda}(\frac{d}{2}Q + x\cdot\nabla Q) + \tilde{\gamma_{s}}Q + i\frac{x_{s}}{\lambda}\cdot\nabla Q + R(\epsilon),
\end{equation}
where $R(\epsilon)$ is formally quadratic in $\epsilon$, and $L = (L_{+}, L_{-})$ is the matrix linearized operator closed to $Q$ which has components:\\

\centerline{$L_{+} = -\Delta + 1 -(1 + \frac{4}{d})Q^{\frac{4}{d}},\quad L_{-} = -\Delta + 1 - Q^{\frac{4}{d}}.$}
\bigskip
A standard approach is to think of equation (\ref{equation en epsilon}) in the following way: it is essentially a linear equation forced by terms depending on the law for the geometrical parameters.\\
Let us observe that the key geometrical parameter is $\lambda$ which measures the size of the solution. Let us then set \\
\centerline{$b = -\frac{\lambda_{s}}{\lambda},$}
\bigskip
and study the simpler version of (\ref{nouvel equation}):\\
\centerline{$iv_{s} + \Delta v - v + v|v|^{\frac{4}{d}} + ib(\frac{d}{2}v + y\cdot\nabla v) = 0.$}
\bigskip

 We look for solutions of the form $v(s,y)=\overline{Q}_{b(s)}(y)$ where the mapping $b\rightarrow Q_{b}$ and the laws for $b(s)$ are the unknown. We take $b$ uniformly small and $Q_{b}|_{b=0}=Q$.
Now injecting $v(s,y)$ into the equation, we get:\\

\centerline{$i\frac{db}{ds}(\frac{\partial{\overline{Q}}_{b}}{\partial b}) + \Delta \overline{Q}_{b(s)} - \overline{Q}_{b(s)} + ib(s)\bigg(\frac{d}{2}\overline{Q}_{b(s)} + y\cdot\nabla \overline{Q}_{b(s)}\bigg) + \overline{Q}_{b(s)}|\overline{Q}_{b(s)}|^{\frac{4}{d}} = 0.$}
\bigskip
We set $\overline{P}_{b(s)} = e^{i\frac{b(s)}{4}}|y|^2\overline{Q}_{b(s)}$ and solve:
\bigskip
\begin{equation}\label{equation pb}
i\frac{db}{ds}(\frac{\partial{\overline{P}}_{b}}{\partial b}) + \Delta \overline{P}_{b(s)} - \overline{P}_{b(s)} +  \bigg(\frac{db}{ds} + b^2(s)\bigg)\frac{|y|^2}{4}\overline{P}_{b(s)}+ \overline{P}_{b(s)}|\overline{P}_{b(s)}|^{\frac{4}{d}} = 0.
\end{equation}
Two remarkable solutions to (\ref{equation pb}) can be obtained as follows:
\begin{itemize}
\item Take $b(s) = 0\quad \text{and}\quad \overline{P}_{b(s)} = Q $, that is the ground state itself.
\item Take $b(s) = b\quad \text{and}\quad \overline{P}_{b(s)} = \overline{P}$ for some non zero constant $b$ and $\overline{P}_{b}$ satisfying:
\begin{equation}\label{eq Pb 1}
\Delta \overline{P}_{b} - \overline{P}_{b} + b^2(s)\frac{|y|^2}{4}\overline{P}_{b}+ \overline{P}_{b}|\overline{P}_{b}|^{\frac{4}{d}} = 0.
\end{equation}
 \end{itemize}
     The solutions to this non linear elliptic equation are those who produce the explicit self similar profiles solutions to this equation:
     \begin{equation}\label{eq en Qb}
      \Delta \overline{Q}_{b} - \overline{Q}_{b}  + \overline{Q}_{b} |\overline{Q}_{b} |^{\frac{4}{d}} + ib(\frac{d}{2}\overline{Q}_{b} + y\cdot\nabla \overline{Q}_{b} ) = 0.
      \end{equation}
      A simple way to see this is to recall that we have set $b =-\frac{\lambda_{s}}{\lambda}$. Hence from $\frac{ds}{dt} = \frac{1}{\lambda^2}$,\\

      \centerline{$ b = -\frac{\lambda_{s}}{\lambda} = -\lambda \lambda_{t}$ ie $\lambda(t) = \sqrt{2b(T-t)},$}
      \bigskip

      which is the scaling law for the blow up speed.

      %Now a crucial point (see \cite{Russell1}) is  that solutions to $(\ref{equation pb})$ never belong to $L^2$ from a logarithmic divergence at infinity:\\

      %\centerline{$|P_{b}(y)| \sim \frac{C(P_{b})}{|y|^{\frac{d}{2}}}$ as $|y| \rightarrow +\infty$}
      %\bigskip
      %This behavior is a consequence of the oscillations induced by the linear group after the turning point $|y| \geq \frac{2}{|b|}$. Nevertheless, in the ball $|y| < \frac{2}{|b|}$, the operator $-\Delta + 1 -\frac{b^2|y|^2}{4}$ is coercive, and no oscillations will take place in this zone.\\

      %Because we track a log-log correction to the self similar law as an upper bound on the blow up speed, the profiles $\overline{Q}_{b} = e^-i\frac{b}{4}|y|^2\overline{P}_{b}$ with $\overline{P}_{b}$ solving (\ref{equation pb}) are natural candidates as refinements of the $Q$ profile in the geometrical decomposition. Nevertheless, as they are not in $L^2$, we need to build a smooth localized version avoiding the non $L^2$ tale, what according to the above discussion is doable in the coercive zone $|y| < \frac{2}{|b|}$.\\
      \section{Setting of the bootstrap}
In this section, we recall some fundamental nonlinear objects central to the description
of the log-log regime. We then set up the
bootstrap argument, in the heart of the proof of Theorem $1.2$. The conditions C.I will be to initialize the bootstrap.\\

Based on Propositions 8 and 9 of \cite{Merle3} , we claim:

      \begin{proposition}There exist universal constants $C > 0$, $\eta^{\star} > 0$ such that the following holds true: for all $0 < \eta < \eta^{\star}$, there exist constants $\nu^{\star}(\eta) > 0$, $b^{\star}(\eta) > 0 $ going to zero as $\eta \rightarrow 0$ such that for all $|b| < b^{\star}(\eta)$, setting \\
      \centerline{$R_{b} = \frac{2}{b}\sqrt{1-\eta}, \quad R_{b}^{-} = \sqrt{1-\eta}R_{b}$,}

      \bigskip
      $B_{R_{b}} = \big\{y \in \mathbb{R}^{d}, |y|\leq R_{b}\big\}$, there exists a unique radial solution $\overline{Q_{b}} \in L^2\big(B(0,R)\big)$ to
\begin{equation}
\begin{cases}
%$$ \left\lbrace \begin{array}{l}
\Delta \overline{Q}_{b} - \overline{Q}_{b}  + \overline{Q}_{b} |\overline{Q}_{b} |^{\frac{4}{d}} + ib(\frac{d}{2}\overline{Q}_{b} + y\cdot\nabla \overline{Q}_{b} ) = 0,\\
\overline{P_{b}}= \overline{Q_{b}}e^i\frac{b|y|^2}{4} > 0 \quad in\quad B_{R_{b}},\\
\overline{Q_{b}(0)} \in \big(Q(0) - \nu^{\star}(\eta), Q(0) + \nu^{\star}(\eta)\big),\overline{Q_{b}}(R_{b}) = 0.
%\end{array} $$
\end{cases}
\end{equation}
Moreover, let $\phi_{b}$ be a smooth radially symmetric cut-off function such that $\phi_{b}(x) = 0$ for $|x| \geq R_{b}$ and $\phi_{b}(x) = 1$ for $|x| \leq R_{b}^{-}, 0 \leq \phi_{b}(x) \leq 1$ and set \\

\centerline{$ {Q}_{b}(r) = \overline{Q_{b}}(r)\phi_{b}(r)$}
\bigskip

then\\

\centerline{${Q}_{b} \rightarrow Q \quad as \quad b \rightarrow 0$}
\bigskip
in $L^2(\mathbb{R}^{d})$, and ${Q}_{b}$ satisfies
\bigskip
\begin{equation}
\Delta {Q}_{b} - {Q}_{b}  + {Q}_{b} |{Q}_{b} |^{\frac{4}{d}} + ib(\frac{d}{2}{Q}_{b} + y\cdot\nabla {Q}_{b})  = - \Psi_{b},
\end{equation}
where $\Psi_{b}=2\nabla\phi_{b}\nabla Q_{b} + Q_{b}(\Delta \phi_{b})+ iQ_{b}y\cdot\nabla\phi_{b} +(\phi_{b}^{1 + \frac{4}{d}} - \phi_{b})Q_{b}\left|Q_{b}\right|^{\frac{4}{d}},$\\
with 
$$supp\big(\Psi_{b}\big) \subset \big\{R_{b}^{-} \leq |y| \leq R_{b} \big\} \quad and \quad | \Psi_{b}|_{C^1} \leq e^{-\frac{c}{|b|}}.$$
\bigskip
Eventually, ${Q}_{b}$ has supercritical mass:
\begin{equation}\label{supercritical mass}
\int | Q_{b}|^2 = \int Q^2 + c_{0}b^2 + o(b^2) \quad as \quad b \rightarrow 0,
\end{equation}
\bigskip
for some universal constant $c_{0} > 0$.\\
\end{proposition}

The meaning of this proposition is that one can build localized ${Q}_{b}$ on the ball $B_{R_{b}}$ which are a smooth function of $b$ and approximate $Q$ in a very strong way as $b \rightarrow 0$. These profiles satisfy the self similar equation up to an exponentially small term $\Psi_{b}$ supported around the turning point $\frac{2}{b}$. The proof of this proposition uses standard variational tools in the setting of non linear elliptic problems, and can be found in \cite{Merle3}.\\

Now one can think of making a formal expansion of $Q_{b}$ in terms of $b$, and the first term is non zero:\\

\centerline{$\frac{\partial{Q}_{b}}{\partial b}|_{b=0} = -\frac{i}{4}|y|^2Q.$}

\bigskip
However, the energy of ${Q}_{b}$ is degenerated in $b$ at all orders:
\begin{equation}\label{energie de Qb}
|E({Q}_{b})| \leq e^{-\frac{c}{|b|}},
\end{equation}
\bigskip
for some universal constant $C > 0$.\\
Now given a well-localized function $f$ , we set:
$$
f_{d}= \frac{d}{2}f + y\cdot\nabla f ~~\text{and}~~ f_{dd}=(f_{d})_{d}.$$

Note that integration by part yields:
$$(f_d,g)_{L^2} = -(g,f_d)_{L^2}.$$

We next introduce the outgoing radiation escaping the soliton core according to the following lemma(see Lemma 15 from \cite{Merle3}):\\
\begin{lemma}\label{Linear outgoing radiation}(Linear outgoing radiation) There exist universal constants $C > 0$ and $\eta^{\star} > 0$ such that $\forall ~ 0 < \eta < \eta^{\star}$, there exists $b^{\star}(\eta) > 0$ such that $\forall |b| < b^{\star}(\eta)$, the following holds true: there exists a unique radial solution $\zeta_{b}$ to
\begin{equation}
\begin{cases}
%$$ \left\lbrace \begin{array}{l}
\Delta\zeta_{b} - \zeta_{b} + ib(\zeta_{b})_{d} = -\Psi_{b}\\
\int |\nabla \zeta_{b}|^2 < +\infty.
%\end{array} $$
\end{cases}
\end{equation}
\bigskip
Moreover, let
\begin{equation}\label{gammab}
\Gamma_{b} = \lim_{|y|\rightarrow +\infty}|y|^d|\zeta_{b}(y)|^2,
\end{equation}
\bigskip
then there holds
\begin{equation}\label{estimation de gammab}
e^{-(1 + c\eta)\frac{\pi}{|b|}} \leq \Gamma_{b}\leq e^{-(1 - c\eta)\frac{\pi}{|b|}}.
\end{equation}
\end{lemma}
\bigskip

We recall that the solution $u(t)$ admits the decomposition :\\

\centerline{$u(t,x)=\frac{1}{\lambda(t)}(Q_{b(t)} + \epsilon)(t,\frac{x-x(t)}{\lambda(t)})e^{i\gamma(t)}$}
\bigskip
where the geometrical parameters are uniquely defined through some orthogonality conditions(see later):
Let us assume the following uniform controls on $[0,T]$:
\begin{itemize}
\item {Control of $b(t)$}
\begin{equation}\label{b petit}
 b(t) > 0 ~ \text{and} ~ b(t) < 10 b(0).
\end{equation}

\item Control of $\lambda$:
\begin{equation}\label{control of lambda}
\lambda(t) \leq e^{-{e^\frac{\pi}{100b(t)}}}
\end{equation}
and the monotonicity of $\lambda$:
\begin{equation}\label{monotonicity}
\lambda(t_2) \leq\frac{3}{2} \lambda(t_{1}), \forall~0\leq t_{1} \leq t_{2} \leq T.
\end{equation}
Let $k_{0} \leq k_{+} $ be an integers and $T^{+} \in [0,T]$  such that
\begin{equation}\label{lambda 0 et lambda T}
\frac{1}{2^{k_{0}}} \leq \lambda(0) \leq \frac{1}{2^{k_{0}-1}}, \frac{1}{2^{k_{+}}} \leq \lambda(T^{+}) \leq \frac{1}{2^{k_{+}-1}}
\end{equation}
and for $k_{0} \leq k \leq k_{+} $, let $t_{k}$ be a time such that
\begin{equation}
\lambda(t_{k}) = \frac{1}{2^k},
\end{equation}
then we assume the control of the doubling time interval:
\begin{equation}\label{tk}
t_{k+1} - t_{k} \leq k \lambda^2(t_{k}).
\end{equation}
\item control of the excess of mass:
\begin{equation}\label{controlepsilon}
\int\left|\nabla \epsilon(t)\right|^2 + \int\left|\epsilon(t)\right|^2e^{-\left|y\right|} \leq \Gamma_{b(t)}^{\frac{1}{4}}.
\end{equation}
\end{itemize}
\bigskip
The following proposition ensures that (\ref{control of lambda})-(\ref{controlepsilon}) determine a trapping region for the flow.
We will prove this proposition in \textbf{section 7} (\textbf{Part 7.3}).
\qed
 \begin{proposition}\label{bootstrap lemma}Assuming that (\ref{b petit})-(\ref{controlepsilon}) hold, then the following controls are also true:
\begin{equation}\label{bestpluspetit}
 b > 0~\text{and}~ b(t) < 5 b(0).
\end{equation}

 \begin{equation}\label{control of lambda 2}
 \lambda(t) \leq e^{-{\frac{\pi}{10b(t)}}}
 \end{equation}
 \begin{equation}\label{monotonicity 2}
 \lambda(t_2) \leq\frac{5}{4} \lambda(t_{1}), \forall \quad0\leq t_{1} \leq t_{2} \leq T
\end{equation}
\begin{equation}\label{tk}
t_{k+1} - t_{k} \leq \sqrt{k} \lambda^2(t_{k})
\end{equation}
\begin{equation}\label{control of epsilon 2}
\int\left|\nabla \epsilon(t)\right|^2 + \int\left|\epsilon(t)\right|^2e^{-\left|y\right|} \leq \Gamma_{b(t)}^{\frac{2}{3}}.
\end{equation}
\end{proposition}
\section{Control of the energy and the kinetic momentum}
We recall the Strichartz estimates. An ordered pair $(q, r )$ is called admissible if $\frac{2}{q} + \frac{d}{r} = \frac{d}{2}$, $2 < q \leq \infty $
 We define the Strichartz norm of functions $u : [0, T]\times\mathbb{R}^{d} \longmapsto \mathbb{C}$ by:
\begin{equation}\label{us0}
\left\|u\right\|_{S^{0}([0,T]\times\mathbb{R}^{d})} = \sup_{(q,r) admissible}\left\|u\right\|_{L^{q}_{t}L^{r}_{x}([0,T]\times\mathbb{R}^{d})}
\end{equation}
and
\begin{equation}\label{us1}
\left\|u\right\|_{S^{1}([0,T]\times\mathbb{R}^{d})} =  \sup_{(q,r) admissible}\left\|\nabla u\right\|_{L^{q}_{t}L^{r}_{x}([0,T]\times\mathbb{R}^{d})}
\end{equation}
We will sometimes abbreviate $S^i([0,T]\times \mathbb{R}^{2})$ with $S^{i}_{T}$  or $S^{i}[0, T]$, $i= 1,2$. Let us denote the H\"older dual exponent of $q$ by $q^{\prime}$ so that $\frac{1}{q} + \frac{1}{q^{\prime}} = 1$. The Strichartz estimates may be
expressed as:
\begin{equation}\label{estimationschwartz}
\displaystyle{\left\|u\right\|_{S^{0}_{T}} \lesssim \left\|u_0\right\|_{L^{2}} + \left\|(i\partial_{t} + \Delta )u\right\|_{L^{q^{\prime}}_{t}L^{r^{\prime}}_{x}}}
\end{equation}
where $(q, r )$ is any admissible pair.
Now we will derive an estimate on the energy, to check that it remains small with respect to $\lambda^{-2}$:

\begin{lemma}\label{control energie}
Assuming that (\ref{control of lambda})-(\ref{controlepsilon}) hold, then the energy and kinetic momentum are controlled on $[0, T^+]$ by:\\
\begin{equation}\label{control de lenergie}
\left|E(u(t))\right|\lesssim \big(\text{log}\big(\lambda(t)\big)\big)^2,
\end{equation}

\begin{equation}\label{control de moment 1}
\left|P(u(t))\right|\leq \left|P(u_{0})\right|.
\end{equation}
\end{lemma}
To prove this lemma, we shall need the following one:\\
\begin{lemma}\label{us1} 
Let $u$ be a solution of $(\ref{NLS})$ emanating for $u_{0}$ in $H^1$. Then $u$ $\in$ $C([0,\Delta T], H^1)$ where $\Delta T = \left\|u_{0}\right\|^{\frac{d-4}{d}}_{L^2}\left\|u_{0}\right\|_{H^{1}}^{-2}$, and we have the following control\\
\centerline{
$\left\|u\right\|_{{S^{0}[t,t+\Delta T]}} \leq 2\left\|u_{0}\right\|_{L^2}$ , \quad $\left\|u\right\|_{{S^{1}[t,t+\Delta T]}} \leq 2 \left\|u_{0}\right\|_{H^1(\mathbb{R}^{d})}.$}
\end{lemma}
\textbf{\underline{proof:}}
For all $v \in S(\mathbb{R}^{d})$ we have:

$$\left\| v\right\|_{S^{0}[t,t+\Delta T]} 
\lesssim \left\|v(0)\right\|_{L^2} + \left\|\left(i\partial_{t} + \Delta\right)v\right\|_{{L^{q^{\prime}}_{\Delta T}}L_{x}^{r^{\prime}}},$$

$\frac{2}{q} + \frac{d}{r} = \frac{d}{2},~ 2 < q < \infty$, $\frac{1}{q} +\frac{1}{q^\prime} = 1.$\\
In particular,\\
$$\left\|\int_{0}^{t} e^{i(t-s)\Delta}|u|^{\frac{4}{d}}u\right\|_{S^{1}_{t}} \lesssim \left\||u|^{\frac{4}{d}}\nabla u\right\|_{{L^1_t}L^2_x}.$$
Using the  H\"older inequality we obtain:
\begin{equation}
\displaystyle{\big(\int\left|u\right|^{\frac{8}{d}}\left|\nabla u\right|^2\big)^{\frac{1}{2}} \leq \big(\int \left|u\right|^{\frac{2(4+d)}{d}}\big)^{\frac{2}{4+d}}\big(\int \left|\nabla u\right|^{\frac{2(4+d)}{d}} \big)^{\frac{d}{2(4+d)}}.}\nonumber
\end{equation}
Integrating in time and applying again H\"older inequality we get:
\begin{align}
\left\|\left|u\right|^{\frac{4}{d}}\nabla u\right\|_{{L^{1}([t,t+\Delta T])}L^{2}(\mathbb{R}^{d})}&\leq \bigg(\int\big(\int \left|u\right|^{\frac{2(4+d)}{d}}dx\big)^{\frac{2}{4+d}\frac{4+d}{4}}dt\bigg)^{\frac{4}{4+d}}\nonumber\\
&\times \bigg(\int\big(\int \left|\nabla u\right|^{\frac{2(4+d)}{d}}dx\big)^{\frac{d}{2(4+d)}\frac{4+d}{d}}dt\bigg)^{\frac{d}{4+d}}.\nonumber
\end{align} 
Thus:
\begin{equation}
\displaystyle{\left\|\left|u\right|^{\frac{4}{d}}\nabla u\right\|_{L^{1}([t,t+\Delta T])L^{2}(\mathbb{R}^{d})} \leq \left\|u\right\|^{\frac{4}{d}}_{L^{\frac{4+d}{d}}([t,t+\Delta T])L^{\frac{8+2d}{d}}(\mathbb{R}^{d})}\left\|\nabla u\right\|_{L^{\frac{4+d}{d}}([t,t+\Delta T])L^{\frac{8+2d}{d}}(\mathbb{R}^{d})}}.\nonumber
\end{equation}
But $(\frac{4+d}{d}, \frac{8+2d}{d})$ is admissible, thus we have:
\begin{equation}
\displaystyle{\left\|\left|u\right|^{\frac{4}{d}}\nabla u\right\|_{{L^{1}([t,t+\Delta T])}L^{2}(\mathbb{R}^{d})} \leq \left\|u\right\|^{\frac{4}{d}}_{L^{\frac{4+d}{d}}([t,t+\Delta T])L^{\frac{8+2d}{d}}(\mathbb{R}^{d})}\left\|u\right\|_{S^{1}[t,t+\Delta T]}}.\nonumber
\end{equation}
By Sobolev we have:
\begin{align}
\left\|u\right\|_{L^{\frac{4+d}{d}[t,t+\Delta T]}L^{\frac{8+2d}{d}}(\mathbb{R}^{d})} &\lesssim \left\|u\right\|_{L^{\frac{4+d}{d}}([t,t+\Delta T])H^{\frac{2d}{d + 4}}(\mathbb{R}^{d})}\nonumber\\
&\leq (\Delta T)^{\frac{d}{d+4}}\left\|u\right\|_{L^{\infty}([t,t+\Delta T])H^{\frac{2d}{d + 4}}(\mathbb{R}^{d})}.\nonumber
\end{align}
Now by interpolation we obtain for $d=1,2,3,4$:
\begin{equation}
\displaystyle{\left\|u\right\|_{L^{\frac{4+d}{d}}([t,t+\Delta T])L^{\frac{8+2d}{d}}(\mathbb{R}^{d})} \leq (\Delta T)^{\frac{d}{d+4}}\left\|u\right\|_{L^{\infty}([t,t+\Delta T])L^{2}(\mathbb{R}^{d})}^{\frac{4-d}{d+4}}\left\|u\right\|_{L^{\infty}([t,t+\Delta T])H^{1}(\mathbb{R}^{d})}^{\frac{2d}{d+4}}.}\nonumber
\end{equation}
But since according to (\ref{mass a}), $\left\|u\right\|_{L^{\infty}([t,t+\Delta T])L^2(\mathbb{R}^{d})} \leq \left\|u_0\right\|_{L^2(\mathbb{R}^{d})}$, we finally get:\\
\centerline{$\left\|u\right\|_{S^{1}[t,t+\Delta T]} \leq \left\|u(t)\right\|_{H^{1}(\mathbb{R}^{d})} + (\Delta T)^{\frac{d}{d+4}}\left\|u_0\right\|_{L^2(\mathbb{R}^{d})}^{\frac{4-d}{d+4}}\left\|u\right\|_{S^{1}_{\Delta T}}^{\frac{2d}{d+4}}\left\|u\right\|_{S^{1}[t,t+\Delta T]},$} 
we deduce that for $\Delta T \leq C \left\|u_0\right\|_{L^2(\mathbb{R}^{d})}^{\frac{d-4}{d}}\left\|u(t)\right\|^{-2}_{H^1(\mathbb{R}^{d})}$:\\

\centerline{$\left\|u\right\|_{S^{1}[t, t + \Delta T]} \leq 2 \left\|u(t)\right\|_{H^1(\mathbb{R}^{d})}.$}

In the same way, $\left\|u\right\|_{S^{0}[t,t+\Delta T]} \lesssim \left\|u_{0}\right\|_{L^2(\mathbb{R}^{d})} + (\Delta T)^{\frac{d}{d+4}}\left\|u_0\right\|^{\frac{4-d}{4+d}}_{L^2(\mathbb{R}^{d})} \left\|u\right\|_{S^{1}[t,t+\Delta T]}^{2\frac{d}{d+4}}\left\|u\right\|_{S^{0}[t,t+\Delta T]}$,\\ but $(\Delta T)^{\frac{d}{d+4}}\left\|u_0\right\|_{L^2(\mathbb{R}^{d})}^{\frac{4-d}{4+d}}\left\|u\right\|_{S^{1}[t,t+\Delta T]}^{2\frac{d}{d+4}} \leq \frac{1}{2}$ we obtain that:\\

\centerline{$\left\|u\right\|_{S^{0}[t,t+\Delta T]} \leq 2 \left\|u_{0}\right\|_{L^2(\mathbb{R}^{d})}~~\text{and}~~\left\|u\right\|_{S^{1}[t,t+\Delta T]} \leq 2 \left\|u_{0}\right\|_{H^1(\mathbb{R}^{d})}.$}
\bigskip
Let us return to the proof of the lemma \ref{control energie}:\\
According to (\ref{tk}) each interval $[t_k, t_{k+1}]$, can be divided into $k$ intervals, $[\tau_{k}^{j},\tau_{k}^{j+1}]$  such that the estimates of the previous lemma are true.
From (\ref{derivee de lenergie}), we thus deduce that:
\begin{align}
\int_{t_{k}}^{t_{k+1}}\frac{d}{dt}E(u(t))&\leq \int_{t_{k}}^{t_{k+1}}\int\left|u\right|^{\frac{4}{d}+2}\nonumber\\
&\leq\sum_{j=1}^{k}\int_{\tau_{k}^{j}}^{\tau_{k}^{j+1}}\int\left|u\right|^{\frac{4}{d}+2}\nonumber\\
&\leq \sum_{j=1}^{k}\left\|u\right\|^{\frac{4}{d}+2}_{{S^0}_{[\tau_{k}^{j},\tau_{k}^{j+1}]}}.\nonumber
\end{align}
Since $(\frac{4}{d}+2,\frac{4}{d}+2)$ is admissible. Using that $\left\|u\right\|_{S^{0}[t,t+\Delta T]} \leq 2 \left\|u_{0}\right\|_{L^2(\mathbb{R}^{d})}$ independantly of $t$, we obtain finally :
\begin{equation}
\displaystyle{\int_{t_{k}}^{t_{k+1}}\frac{d}{dt}E(u(t)) \lesssim k.}\nonumber
\end{equation}
Summing from $k_{0}$ to $k_{+}$, we obtain 
\begin{equation}
\displaystyle{\int_{0}^{T^+}\frac{d}{dt}E(u(t)) \lesssim \left|\text{log}(\lambda(T^+))\right|^2.}\nonumber
\end{equation}
But $E(u(T^+)) = E(u(0)) + \displaystyle{\int}_{0}^{t}\frac{d}{ds}E(u(s))ds$,
 since  $\left|E(u(0))\right| \leq \frac{1}{\sqrt{\lambda(0)}}$ then %$\lambda(t) \leq \frac{3}{2} \lambda(0)$
 we obtain $\left|E(u(T^+))\right| \lesssim \left|\text{log}(\lambda(T^+))\right|^2$. This completes the proof of (\ref{control de lenergie}). \\
Now (\ref{control de moment 1}) follow directly from (\ref{moment}). For sake of completness, let us prove (\ref{moment}):\\
Suppose first that $u(t)$ is very regular (for example $D(\mathbb{R}^{d})$)
\begin{align}
\frac{d}{dt}P(u(t)) &= Im \bigg(\int\overline{u_{t}}\nabla{u}dx + \int{\overline{u}\nabla{u_{t}}}\bigg)\nonumber\\
&=Im \bigg(\int\overline{u_{t}}\nabla{u}dx - \int u_{t}\nabla{\overline{u}}\bigg)\nonumber\\
&=Im(2iIm\int{\overline{u_{t}}\nabla u})\nonumber\\
&=2Im \int{\overline{u_{t}}\nabla udx}= 2Im(i\int{\overline{iu_{t}}\nabla u}) = 2Re\int\overline{iu_{t}}\nabla u\nonumber\\
&=2Re\bigg(\int(-\Delta\overline{u} - \left|u\right|^{\frac{4}{d}}\overline{u} + ia\overline{u})\nabla udx\nonumber\\
&=-2Re(\int\Delta\overline{u}\nabla u) -2Re\int{\left|u\right|^{\frac{4}{d}}\overline{u}\nabla u} + 2aRe\int{i\overline{u}\nabla u}.\nonumber
\end{align}
It is easy to prove that: $-2Re(\displaystyle{\int}\Delta\overline{u}\nabla u) = \displaystyle{\int}\nabla \left|\nabla u\right|^2$, and $-2Re\displaystyle{\int}{\left|u\right|^{\frac{4}{d}}\overline{u}\nabla u} = -\frac{d}{d+2}\displaystyle{\int}\nabla(\left|u\right|^{\frac{4}{d} + 2}).$ \\
But $2aRe\displaystyle{\int}i\overline{u}\nabla{u} = -2aIm\displaystyle{\int}\overline{u}\nabla u = -2aP(u(t))$ we obtain:\\
$\frac{d}{dt}P(u(t)) =\displaystyle{\int}\nabla (\left|\nabla u\right|^2 - \frac{d}{d+2}\left|u\right|^{\frac{4}{d} + 2}) - 2aP(u(t))$, But $\displaystyle{\int}\bigg(\nabla (\left|\nabla u\right|^2 - \frac{d}{d+2}\left|u\right|^{\frac{4}{d} + 2})\bigg) = 0$, we obtain finally:
$$\frac{d}{dt}P(u(t))= -2a p(u(t))~~\text{and}~~P(u(t)) = e^{-2at}P(u_{0}).$$
Now we take $u_{0}$ in $H^1(\mathbb{R}^{d})$, $u_{0}$ is the limit of a sequence $(u_{0n})$ in $D(\mathbb{R}^{d})$, for each $u_{0n}$ we denote by $u_{n}$ the solution of $(\ref{NLSa})$ such that $u_{n}(0) = u_{0n}$, we have $P(u_{n}(t)) = e^{-2at}P(u_{0n})$, but $u_{0}\mapsto u$ is continuous from $H^1(\mathbb{R}^{d})$ to $C\big([0,T], H^{1}(\mathbb{R}^{d})\big)$ by passing to limit we obtain $P(u(t)) = e^{-2at}P(u(0))$.\hfill$\Box$

 \section{Booting the log-log regime} 
 Now we are going to prove the Lemma \ref{bootstrap lemma}:
 \\
 First of all we are going to prove the smallness of the $L^2$ norm of $\epsilon(t)$:
 \begin{lemma}
 There exist $\alpha_{0} \ll 1$, such that $\forall t \in [0,T],~ \left\|\epsilon(t)\right\|_{L^2(\mathbb{R}^{d})} < \alpha_{0}$.
 \end{lemma}
 \textbf{Proof:} From (\ref{small norm}) we have $\left\|u_{0}\right\|_{L^2} < \left\|Q\right\|_{L^2} + \gamma_{0}$, with $\gamma_{0}$ very small. By (\ref{mass a})
 \begin{align}
 \left\|u_{0}\right\|_{L^2(\mathbb{R}^{d})}^2 &\geq \left\|u(t)\right\|_{L^2(\mathbb{R}^{d})}^2 \nonumber\\
 &= \left\|Q_{b} + \epsilon\right\|_{L^2(\mathbb{R}^{d})}^2 \nonumber\\
 &= \left\|Q_{b} + \epsilon\right\|_{L^2(B(0,R))}^2 + \left\|Q_{b} + \epsilon\right\|_{L^2(\mathbb{R}^{d}\backslash B(0,R))}^2.\nonumber\\
 \left\|Q_{b} + \epsilon\right\|_{L^2(B(0,R))} &\geq  \left\|Q_{b}\right\|_{L^2(B(0,R))} - \left\|\epsilon\right\|_{L^2(B(0,R))}\nonumber\\
 & \geq \left\|Q\right\|_{L^2(\mathbb{R}^{d})} -  \left\|Q_{b} - Q\right\|_{L^2(\mathbb{R}^{d})}\nonumber\\ 
 & - \left\|Q_{b}\right\|_{L^{2}(\mathbb{R}^{d}\backslash B(0,R))} - \left\|\epsilon\right\|_{L^2(B(0,R))}.\nonumber
 \end{align} 
 From (\ref{controlepsilon}),  we have $\left\|\epsilon\right\|_{L^2(B(0,R))} < \beta$, where  $\beta$ is very small. Moreover $Q_{b}\rightarrow Q$ in $L^2(\mathbb{R}^{d})$, and $\left\|Q_{b}\right\|_{L^2(\mathbb{R}^{d}\backslash B(0,R))} \leq \left\|Q_{b} - Q\right\|_{L^2(\mathbb{R}^{d})} + \left\|Q\right\|_{{L^2(\mathbb{R}^{d}\backslash B(0,R))}}$, where $\left\|Q\right\|_{{L^2(\mathbb{R}^{d}\backslash B(0,R))}} \rightarrow 0$ as $R\rightarrow \infty$.\\
 Therefore: 
 \begin{align}
 \left\|Q_{b} + \epsilon\right\|_{L^2(B(0,R))}\geq \left\|u_{0}\right\|_{L^2(\mathbb{R}^{d})} - \gamma_{0} - \beta - \delta,\nonumber
\end{align}
and\\
\begin{align}
 \left\|Q_{b} + \epsilon\right\|_{L^2(\mathbb{R}^{d}\backslash B(0,R))} &\geq \left\|\epsilon\right\|_{L^2(\mathbb{R}^{d}\backslash B(0,R))} - \left\|Q_{b}\right\|_{L^2(\mathbb{R}^{d}\backslash B(0,R))}\nonumber\\
  &\geq  \left\|\epsilon\right\|_{L^2(\mathbb{R}^{d}\backslash B(0,R))} - \delta,\nonumber
  \end{align}
   where $\delta \rightarrow 0$ as $b\rightarrow 0$ and $R \rightarrow \infty$. We obtain finally:\\
   \centerline{ $\left\|u_{0}\right\|_{L^2(\mathbb{R}^{d})}^2 \geq  \left\|\epsilon\right\|_{L^2(\mathbb{R}^{d}\backslash B(0,R))}^2 + \left\|u_{0}\right\|_{L^2(\mathbb{R}^{d})}^2 - \alpha_{0}^2$,}\\
    where $ \alpha_{0}  \rightarrow 0$ as $\gamma_0 \rightarrow 0$.\\
    This completes the proof.\hfill$\Box$ 
 \subsection{Control of the geometrical parameters} 
 Let us now write down the equation satisfied by $\epsilon$ in rescaled variables. To simplify the notations, we note \\
\centerline{${Q}_{b} = \Sigma + i\Theta ~~ ,~~ \epsilon = \epsilon_1 + i\epsilon_2 ~~\text{and}~~ \Psi_{b}=Re(\Psi) + iIm(\Psi)$,}\\
in terms of real and imaginary parts.\\
 We have: $\forall s \in \mathbb{R}^{+}, \forall y \in \mathbb{R}^{d}$,
\begin{equation}\label{eq 1 de epsilon}\begin{array}[t]{lll}
b_{s}\frac{\partial\Sigma}{\partial b} + \partial_{s}\epsilon_{1} - M_{-}(\epsilon) + b(\epsilon_{1})_{d}&\!\!\!=\!\!\!&\bigg(\frac{\lambda_{s}}{\lambda} + b\bigg) \Sigma_{d} + \tilde{\gamma}_{s}\Theta + \frac{x_{s}}{\lambda}\cdot\nabla \Sigma\medskip\\ &\!\!\!+\!\!\!&\bigg(\frac{\lambda_{s}}{\lambda} + b\bigg) (\epsilon_{1})_{d} + \tilde{\gamma}_{s}\epsilon_{2} + \frac{x_{s}}{\lambda}\cdot\nabla \epsilon_{1}\medskip\\
&\!\!\!+\!\!\!& Im(\Psi) - R_{2}(\epsilon).
\end{array}
\end{equation}
\begin{equation}\label{eq 2 de epsilon}\begin{array}[t]{lll}
b_{s}\frac{\partial \Theta}{\partial b} + \partial_{s}\epsilon_{2} + M_{+}(\epsilon) + b(\epsilon_{2})_{d} &\!\!\!=\!\!\!& \bigg(\frac{\lambda_{s}}{\lambda} + b\bigg){\Theta}_{d} - \tilde{\gamma}_{s}\Sigma + \frac{x_{s}}{\lambda}\cdot\Theta\medskip\\
&\!\!\!+\!\!\!&\bigg(\frac{\lambda_{s}}{\lambda} + b\bigg)(\epsilon_{2})_{d} - \tilde{\gamma}_{s}\epsilon_{1} + \frac{x_{s}}{\lambda}\cdot\nabla \epsilon_{2}\medskip\\
&\!\!\!-\!\!\!& Re(\Psi) + R_{1}(\epsilon).
\end{array}
\end{equation}
\bigskip
With $\tilde{\gamma}_{s}(s) = -1 - \gamma_{s}(s) -ia\lambda^2$. The linear operator close to $Q_{b}$ is now a deformation of the linear operator $L$ close to $Q$ and is $M = (M_{+},M_{-})$ with \\
\\
$M_{+}(\epsilon) = -\Delta\epsilon_{1} + \epsilon_{1} - \bigg(\frac{4\Sigma^2}{d|{Q}_{b}|^2} + 1 \bigg)|Q_{b}|^{\frac{4}{d}}\epsilon_{1} - \bigg( \frac{4\Sigma\Theta}{d|{Q}_{b}|^2}|{Q}_{b}|^{\frac{4}{d}}\bigg)\epsilon_{2}$,\\
$M_{-}(\epsilon) = -\Delta\epsilon_{2} + \epsilon_{2} - \bigg(\frac{4\Sigma^2}{d|\tilde{Q}_{b}|^2} + 1 \bigg)| Q_{b}|^{\frac{4}{d}}\epsilon_{2} - \bigg( \frac{4\Sigma\Theta}{d|{Q}_{b}|^2}|{Q}_{b}|^{\frac{4}{d}}\bigg)\epsilon_{1}$.\\
The formally quadratic in $\epsilon$ interaction terms are:\\
$R_{1}(\epsilon)= (\epsilon_{1} + \Sigma)|\epsilon + {Q}_{b}|^{\frac{4}{d}} - \Sigma| Q_{b}|^{\frac{4}{d}}- \bigg(\frac{4\Sigma^2}{d|{Q}_{b}|^2} + 1 \bigg)| Q_{b}|^{\frac{4}{d}}\epsilon_{1} - \bigg( \frac{4\Sigma\Theta}{d|{Q}_{b}|^2}|{Q}_{b}|^{\frac{4}{d}}\bigg)\epsilon_{2}.$\\
$R_{2}(\epsilon)=(\epsilon_{2} + \Sigma)|\epsilon + {Q}_{b}|^{\frac{4}{d}} - \Sigma|Q_{b}|^{\frac{4}{d}}- \bigg(\frac{4\Sigma^2}{d|{Q}_{b}|^2} + 1 \bigg)| Q_{b}|^{\frac{4}{d}}\epsilon_{2} - \bigg( \frac{4\Sigma\Theta}{d|{Q}_{b}|^2}|{Q}_{b}|^{\frac{4}{d}}\bigg)\epsilon_{1}.$\\
We note $s(0)$ by $s_0$ and $s(T^+)$ by $s^+$, now we have the following lemma:
\begin{lemma}\label{control of the geometrical}(Control of the geometrical parameters)
For all $ s \in [s_{0},s^{+}]$, there holds:
\begin{itemize}
\item Estimates induced by the control of energy and momentum:
\begin{align}\label{control 1}
&|2\left(\epsilon_{1}, \Sigma + b\Theta_{d} - Re(\Psi_{b})) + 2(\epsilon_{2}, \Theta - b\Sigma_{d}- Im(\Psi_{b})\right)|\nonumber\\
&\leq\delta_{0}(\int{\left|\nabla \epsilon\right|^2} +\int \left|\epsilon\right|^{2}e^{-\left|y\right|})^{\frac{1}{2}} + \Gamma_{b(s)}^{1-c\eta} + \lambda^{2}\left|E(u(t))\right|.
 \end{align}
\begin{equation}\label{control 2}
\left|(\epsilon_{2},\nabla Q)\right| \leq \delta_{0}\left\|\nabla \epsilon(s)\right\|_{L^2(\mathbb{R}^{d})} + \lambda\left|P(u)\right|.
\end{equation}
\item Estimates on the modulation parameters:
\begin{equation}\label{control 3}
\left|\frac{\lambda_{s}}{\lambda} + b\right| + \left|b_{s}\right| \leq C \left(\int{\left|\nabla\epsilon(s)\right|^2} + \int{\left|\epsilon(s)\right|^2e^{-\left|y\right|}}\right)^{\frac{1}{2}} + \Gamma_{b(s)}^{1-c\eta}.
\end{equation}
\begin{equation}\label{control 4}
\left|\tilde{\gamma_{s}} - \frac{(\epsilon_{1},L_{+}Q_{dd})}{\left\| Q_{d}\right\|_{L^2}^{2}}\right| + \left|\frac{x_{s}}{\lambda}\right| \leq \delta_{0}\left(\int{\left|\epsilon(s)\right|^{2}e^{-\left|y\right|}}\right)^{\frac{1}{2}} + \Gamma_{b(s)}^{1-c\eta}.
\end{equation}
Here $\delta_{0}$ is a small constant $\delta_0 \ll 1$.
\end{itemize}
\end{lemma}
We will need the following lemma ( for the proof see \cite{Merle2} ).
\begin{lemma}\label{nouvellemme}(Control of nonlinear interactions). Let $P(y)$ a polynomial and
integers $0 \leq k \leq 3$, $0 \leq l \leq 1$, $0 \leq m \leq 2$, then for some function
 $\delta(\alpha_{0})\longrightarrow 0$ as $\alpha_{0} \longrightarrow 0$,
\begin{itemize}
\item$\left|\big(\epsilon,P(y)\frac{d^k}{dy^k}Q_{b}(y)\big)\right| \leq C_{P,k}(\displaystyle{\int}{\left|\epsilon(s)\right|^{2}e^{-\left|y\right|}})^{\frac{1}{2}},$
\item$\left|\big(\epsilon,P(y)\frac{d^k}{dy^k}(Q_{b}(y)- Q(y)\big)\right| \leq \delta(\alpha_{0})(\displaystyle{\int}{\left|\epsilon(s)\right|^{2}e^{-\left|y\right|}})^{\frac{1}{2}},$
\item$\displaystyle{\int}\left|\epsilon\right|\left|P(y)\frac{d^m}{dy^m}\frac{\partial{Q_{b}}}{\partial b}\right| \leq (\displaystyle{\int}{\left|\epsilon(s)\right|^{2}e^{-\left|y\right|}} + \displaystyle{\int}\left|\nabla \epsilon\right|^2dy)^{\frac{1}{2}},$

\item$\left|\big(R(\epsilon),P(y)\frac{d^k}{dy^k}Q_{b}(y)\big)\right| \leq C(\int{\left|\epsilon(s)\right|^{2}e^{-\left|y\right|}dy} + \displaystyle{\int}\left|\nabla \epsilon\right|^2dy)^{\frac{1}{2}},$
\item$\displaystyle{\int}\left|F(\epsilon)\right| + \left|(\tilde{R_{1}}(\epsilon), \Sigma_{d})\right| +\left|(\tilde{R_{2}}(\epsilon), \Theta_{d})\right| \leq \delta(\alpha_{0})(\displaystyle{\int}{\left|\epsilon(s)\right|^{2}e^{-\left|y\right|}dy} + \displaystyle{\int}\left|\nabla \epsilon\right|^2dy)^{\frac{1}{2}},$
\item$\left|\big(P(y)\frac{d^l}{dy^l}\Psi,\frac{d^k}{dy^k}Q_{b}(y)\big)\right| + \left|\big(\epsilon, P(y)\frac{d^l}{dy^l}\Psi\big)\right| \leq e^{-\frac{C}{|b|}},
$
\item$\left|\big(\frac{\partial Q_{b}}{\partial b},P(y)\frac{d^k}{dy^k}Q_{b}(y)\big) + \big(i\frac{|y|^2}{4}Q,P(y)\frac{d^k}{dy^k}Q(y)\big)\right|\leq \delta(\alpha_{0}).
$
\end{itemize}
\end{lemma} 
\textbf{Proof of Lemma \ref{control of the geometrical}:}
To prove (\ref{control 1}), we rewrite the expression of energy in the $\epsilon$ variable ($\epsilon = e^{i\gamma(t)}\lambda^{\frac{d}{2}}(t)u\big(t, \lambda(t)x + x(t)\big) - Q_b$):
\begin{align}\label{demon control 1}
&2(\epsilon_{1}, \Sigma + b\Theta_{d} - Re(\Psi)) + 2(\epsilon_{2}, \Theta - b\Sigma_{d} - Im(\Psi))\nonumber\\
&= 2E(Q_{b}) - 2\lambda^2E(u(t))\nonumber\\
&+\int{\left|\nabla\epsilon\right|^2} - \int(\frac{4\Sigma^2}{d\left|Q_{b}\right|^2} + 1)\left|Q_{b}\right|^{\frac{4}{d}}\epsilon_{1}^2\nonumber\\
&-\int(\frac{4\Theta^2}{d\left|Q_{b}\right|^2} + 1)\left|Q_{b}\right|^{\frac{4}{d}}\epsilon_{2}^2\nonumber\\
&-8\int\frac{\Sigma\Theta}{d\left|Q_{b}\right|^2}\left|Q_{b}\right|^{\frac{4}{d}}\epsilon_{1}\epsilon_{2} -\frac{2}{2+\frac{4}{d}}\int F(\epsilon).
\end{align}
With
\begin{align}
F(\epsilon) &= \left|\epsilon + Q_{b}\right|^{\frac{4}{d}+2} -\left|Q_{b}\right|^{\frac{4}{d} + 2} - \bigg(\frac{4}{d} + 2\bigg)\frac{\left|Q_{b}\right|^{\frac{4}{d} + 2}}{\left|Q_{b}\right|^2}(\Sigma\epsilon_{1} + \Theta\epsilon_{2})\nonumber\\
&-\epsilon_{1}^{2}\frac{\left|Q_{b}\right|^{\frac{4}{d} + 2}}{\left|Q_{b}\right|^4}\bigg(\big(\frac{2}{d}+1\big)\big(\frac{4}{d}+1\big)\Sigma^2 + \big(\frac{2}{d}+1\big)\Theta^2\bigg)\nonumber\\
&-\epsilon_{2}^{2}\frac{\left|Q_{b}\right|^{\frac{4}{d} + 2}}{\left|Q_{b}\right|^4}\bigg(\big(\frac{2}{d}+1\big)\big(\frac{4}{d}+1\big)\Theta^2 + \big(\frac{2}{d}+1\big)\Sigma^2\bigg)\nonumber\\
&-\epsilon_{1}\epsilon_{2}\frac{\left|Q_{b}\right|^{\frac{4}{d} + 2}}{\left|Q_{b}\right|^4}\frac{8}{d}\bigg(\frac{2}{d} + 1\bigg)\Sigma\Theta.\nonumber
\end{align}
 (\ref{control 1}) then follows from Lemma \ref{nouvellemme} \big(we estimate the terms in  (\ref{demon control 1}) using Lemma \ref{nouvellemme}, and we obtain (\ref{control 1})\big).\\
Now to prove (\ref{control 2}), we rewrite the expression of the moment in the $\epsilon$ variable:
\begin{align}
P(u(t)) = Im\int\big({\nabla u\overline{u}}\big)&=\frac{1}{\lambda}Im\bigg(\int(\nabla \epsilon + \nabla Q_{b})\overline{(\epsilon + Q_{b})}\bigg)\nonumber\\
&=\frac{1}{\lambda}\bigg(Im\big(\int\nabla \epsilon\overline{\epsilon}\big) - 2(\epsilon_{2},\nabla\Sigma) + 2(\epsilon_{1}, \nabla\Theta)\bigg),\nonumber
\end{align}
so that\\
$$2(\epsilon_{2}, \nabla \Sigma)=2(\epsilon_{1}, \nabla \Theta) + Im\big(\int{\nabla \epsilon\overline{\epsilon}}\big) -\lambda P(u(t)).$$
From $\Theta_{b=0}=0$ and the smallness of the $L^2$ norm of $\epsilon(t)$ and the control(\ref{moment}) of the momentum  we obtain (\ref{control 2}).\\   
The prove (\ref{control 3})-(\ref{control 4}), it suffices to follow the proof of Lemma 3 in \cite{Merle2}.\qed

\bigskip
Now let
\begin{align}\label{R1tilde}
 \tilde{R_{1}}(\epsilon) &= R_{1}(\epsilon) -\epsilon_{1}^2\frac{\left|Q_{b}\right|^{\frac{4}{d}}}{\left|Q_{b}\right|^4}\bigg(\frac{2}{d}(\frac{4}{d} +1)\Sigma^3 + \frac{6}{d}\Sigma\Theta^2\bigg)\nonumber\\
 &-\epsilon_{2}^2\frac{\left|Q_{b}\right|^{\frac{4}{d}}}{\left|Q_{b}\right|^4}\bigg(\frac{2}{d}\Sigma^3 + \frac{2}{d}(\frac{4}{d} - 1)\Sigma\Theta^2\bigg)\nonumber\\
 &-\frac{4}{d}\frac{\left|Q_{b}\right|^{\frac{4}{d}}}{\left|Q_{b}\right|^4}\epsilon_{1}\epsilon_{2}\bigg((\frac{4}{d}-1)\Sigma^2\Theta + \Theta^3\bigg),
 \end{align}
 \begin{align}\label{R2tilde}
 \tilde{R_{2}}(\epsilon) &= R_{2}(\epsilon) -\epsilon_{2}^2\frac{\left|Q_{b}\right|^{\frac{4}{d}}}{\left|Q_{b}\right|^4}\bigg(\frac{2}{d}(\frac{4}{d} +1)\Theta^3 + \frac{6}{d}\Sigma^2\Theta\bigg)\nonumber\\
 &-\epsilon_{1}^2\frac{\left|Q_{b}\right|^{\frac{4}{d}}}{\left|Q_{b}\right|^4}\bigg(\frac{2}{d}\Theta^3 + \frac{2}{d}(\frac{4}{d} - 1)\Sigma^2\Theta\bigg)\nonumber\\
 &-\frac{4}{d}\frac{\left|Q_{b}\right|^{\frac{4}{d}}}{\left|Q_{b}\right|^4}\epsilon_{1}\epsilon_{2}\bigg((\frac{4}{d}-1)\Sigma\Theta^2 + \Sigma^3\bigg).
 \end{align} 
 We define the two real Shr\"odinger operators:
$$L_1 = -\Delta + \frac{2}{d}(\frac{4}{d}+1)Q^{\frac{4}{d}-1}y\cdot\nabla Q, ~~ L_2 = -\Delta + \frac{2}{d}Q^{\frac{4}{d}-1}y\cdot\nabla Q.$$
 To show the explosion, we will need to control $b_{s}$. Note that our continuous functions $(\lambda, \gamma, x(t), b)$ such that:
$$
\epsilon = e^{i\gamma(t)}\lambda^{\frac{d}{2}}(t)u\big(t, \lambda(t)x + x(t)\big) - Q_b
$$
satisfy the following conditions of orthogonality:
\begin{equation}\label{orth1}
(\epsilon_{1}(t),\Sigma_{d}) + (\epsilon_{2}(t),\Theta_{d}) = 0,
\end{equation}
\begin{equation}\label{orth2}
(\epsilon_{1}(t),y\Sigma) + (\epsilon_{2}(t),y\Theta) = 0,
\end{equation}
\begin{equation}\label{orth3}
-(\epsilon_{1}(t),\Theta_{dd}) + (\epsilon_{2}(t),\Sigma_{dd}) = 0,
\end{equation}
\begin{equation}\label{orth4}
-(\epsilon_{1}(t),\Theta_{d}) + (\epsilon_{2}(t),\Sigma_{d}) = 0.
\end{equation}
For the proof of these conditions see Lemma 2 in \cite{Merle2}, the proof is based  on the implicit function theorem using that $(Q_b)_{b=0} = Q $ \linebreak and $(\frac{\partial Q_b}{\partial b})_{b=0} = -i\frac{\left|y\right|^2}{4}Q$.\\
Now we have the following one:
\begin{proposition}There exist $\delta_{0} > 0$, $C > 0$ and $ 0 < \beta < 2 $ such that:
\begin{equation}\label{control de bs}
b_{s} \geq \delta_{0}\big(\int\left|\nabla \epsilon\right|^2 + \int\left|\epsilon\right|^2e^{-\left|y\right|}\big) - e^{-\frac{C}{b}} - \lambda^{\beta}(s).
\end{equation}
\end{proposition}
\textbf{Proof:}\,To prove this proposition, it suffices to follow the proof of Proposition 3 in \cite{Merle2} and used the control of the energy.
 .$\hfill\Box$\\\\

We will need to refine $Q_{b}$ because $Q_{b}$ is not an exact sel-similar solution. The basic idea is that the profile $Q_{b} + \zeta_{b}$ should be a better approximation of the solution. 
 The problem is now that $\zeta_{b}$ is indeed in $\dot{H}^ 1$, but not in $L^2$, and we  then are not able to estimate the main interaction terms.
We therefore introduce a cut version of the radiation: leave a radial cutoff function: $\chi_{A}(r)=\chi(\frac{r}{A})$ with $\chi(r) = 1$ for $0\leq r \leq 1$ and $\chi(r)=0$ for $r\geq2$.\\
The choice of the parameter $A(t)$ is a crucial issue in our analysis, and is roughly based on two contraints: we want $A$ to be large in order first to enter the radiative zone, i.e., $ \frac{2}{b} \ll A$, and to ensure the slowest possible variations of the $L^2$-norm in the zone $\left|y\right| \geq A$. But we also want $A$ not too large, in particular to keep a good control over local $L^2$-terms of the form $\int_{\left|y\right| \leq A}\left|\epsilon\right|^2$.\\ 
A choice which balances these two contraints is:\\
\centerline{$A = A(t) = e^{2l\frac{1}{b(t)}}$ so that $\Gamma_{b}^{-\frac{l}{2}} \leq A \leq \Gamma_{b}^{-\frac{3l}{2}}$,}
for some parameter $ l > 0 $ small enough to be chosen later and which depends on $\eta$. 
Now let\\ 
\centerline{$\tilde{\zeta} = \chi(\frac{r}{A})\zeta_{b}$,}
Observe that $\tilde{\zeta}$ is now a small Schwartz function thanks to the $A$ localization. we next consider the new variable 
\begin{equation}\label{nouvelle variable}
\tilde{\epsilon} = \epsilon -\tilde{\zeta},
\end{equation}
$\tilde{\zeta}_{b}$ still satisfies the size estimates of Lemma $2.1$ and is moreover in $L^2$ with an estimate 
\begin{equation}\label{estimation su zeta tilde}
\int\left|\tilde{\zeta}\right|^2 \leq \Gamma_{b}^{1-C\eta}.
\end{equation}
The equation satisfied by $\tilde{\zeta}$ is now\\ 
\centerline{$\Delta \tilde{\zeta} - \tilde{\zeta} + ib(\tilde{\zeta})_{d} = \Psi_{b} + F$}
with
\begin{equation}\label{Expression de F}
F=(\Delta\chi_{A})\zeta_{b} + 2\nabla\chi_{A}\cdot\nabla\zeta_{b} +iby\cdot\nabla\chi_{A}\zeta_{b}.
\end{equation}
Now we have the following lemma (see Lemma $4.4$ in \cite{Colliander} for further details):
\begin{lemma}(Virial dispersion in the radiative regime) There exist a constants $\delta_{1} > 0 $, $C > 0$ and $\alpha > 0$ such that:
\begin{align}\label{f1s}
\big(f_{1}(s)\big)_{s} &\geq \delta_{1}\bigg(\int \left|\nabla \tilde{\epsilon}\right|^2 + \int\left|\epsilon(s)\right|^2e^{-\left|y\right|} + \Gamma_{b}\bigg)\nonumber\\
& -\frac{1}{\delta_{1}}\int_{A}^{2A}\left|\epsilon\right|^2 - C\lambda^2E(u(t)).
\end{align}
  with 
  \begin{equation}\label{f1}
  f_{1}(s) = \frac{b}{4}\left\|yQ_{b}\right\|_{L^2}^2 + \frac{1}{2}Im \big(\int(y\cdot\nabla\tilde{\epsilon})\overline{\tilde{\zeta}}\big) + (\epsilon_{2}, \Lambda\tilde{\zeta}_{re}) - (\epsilon_{1}, \Lambda\tilde{\zeta}_{im}).
  \end{equation}
  \end{lemma}
  We now need to control the term $\displaystyle{\int}_{A}^{2A}\left|\epsilon\right|^2$ in (\ref{f1s}). This is achieved by computing the flux of $L^2$-norm escaping the radiative zone. We introduce a radial nonnegative cut off function $\phi(r)$ such that $\phi(r) = 0$ for $r \leq \frac{1}{2}$, $\phi(r) = 1$ for $r \geq 3$, $\frac{1}{4} \leq \phi^{\prime}(r) \leq \frac{1}{2}$ for $1 \leq r \leq 2$, $\phi^{\prime}(r) \geq 0$. We then set \\
  \centerline{$\phi_{A}(s,r) = \phi(\frac{r}{A(s)})$.}
  Moreover, we restrict the freedom on the choice of the parameters $(\eta,l)$ by assuming $ l > C\eta$.
  We have:\\
  \begin{align}
  \begin{cases}
%$$ \left\lbrace \begin{array}{l}
\phi_{A}(r) = 0 \quad \text{for}\quad 0 \leq r \leq \frac{A}{2},\\
\frac{1}{4A} \leq \phi_{A}^{\prime}(r) \leq \frac{1}{2A}\quad \text{for}\quad A \leq r \leq 2A,\\
\phi_{A}(r) = 1 \quad\text{for}\quad r \geq 3A\\
\phi_{A}^{\prime}(r) \geq 0, 0 \leq \phi_{A}(r) \leq 1.
%\end{array} $$
\end{cases}\nonumber
\end{align}
We now claim the following dispersive control at infinity in space ( see Lemma 7 in \cite{Merle4} for the proof):
\begin{lemma}($L^2$ dispersion at infinity in space). For some universal constant $C > 0$, if $s$ large enough:
\begin{equation}\label{control de phiA}
\big(\int \phi_{A}\left|\epsilon\right|^2\big)_{s} \geq \frac{b}{400}\int_{A}^{2A}\left|\epsilon\right|^2 - \Gamma_{b}^{1+Cl} - \Gamma_{b}^{\frac{l}{2}}\int\left|\nabla \epsilon\right|^2 - \frac{\lambda^{2}}{b^2}E(u(t)).
\end{equation}
\end{lemma}

 Note that $\lambda \leq e^{-e^{\frac{C}{b}}}$ thus $\frac{\lambda^2}{b^2}E(u(t)) \leq \lambda^\beta$ with $0 < \beta <  2$ close to 2, thus the last term in this estimation is small with to respect $\lambda $.

\subsection{$L^2$-dispersive constraint on the solution.}
\bigskip
 In this subsection, we derive the
dispersive estimate needed for the proof of the blowup. The virial estimate (\ref{f1s})
corresponds to nonlinear interactions on compact sets. The $ L^2 $linear estimate (\ref{control de phiA})
measures the interactions with the linear dynamic at infinity. We now couple these
two facts through the smallness of the $L^2$-norm, which is a global information in
space.\\
\begin{proposition}For some universal constant $C > 0$ and for $s \geq 0$, the following holds:
\begin{equation}\label{derniere estimation}
(\Im)_{s} \leq -Cb\bigg(\Gamma_{b} + \int\left|\nabla \tilde{\epsilon}\right|^2 +\int\left|\tilde{\epsilon}\right|^2e^{-\left|y\right|} + \int_{A}^{2A}\left|\epsilon\right|^2\bigg) + C\frac{\lambda^2}{b^2}E(u(t)).
\end{equation}
\end{proposition}
with:
\begin{align}\label{f}
\Im(s)&= \big(\int\left|Q_{b}\right|^2 - \int\left|Q\right|^2\big) + 2(\epsilon_{1}, \Sigma) + 2(\epsilon_{2}, \Theta) +\int(1-\phi_{A})\left|\epsilon\right|^2\nonumber\\
&-\frac{\delta_1}{800}\bigg(b\tilde{f_{1}}(b) - \int_{0}^{b}\tilde{f_{1}}(v)dv + b\big((\epsilon_{2},\Lambda\tilde{\zeta_{re}}) - (\epsilon_{1}, \Lambda\tilde{\zeta_{im}})\big)\bigg).
\end{align}
Here $c>0$ denotes some small enough universal constant and:
\begin{equation}\label{f1tilde}
\tilde{f_{1}}(b)= \frac{b}{4}\left|yQ_{b}\right|_{2}^2 + \frac{1}{2}Im\big(\int(y\cdot\nabla\tilde\zeta)\overline{\tilde{\zeta}}\big).
\end{equation}
\begin{remark}
Here the range of parameters is more restricted and yields: there exist
$\eta^{\ast}$, $l^{\ast}$, $C_{0} > 0$ such that $\forall ~ 0 < \eta < \eta^{\ast}$, $\forall~ 0 < l < l^{\ast}$ such that $l > C_{0}\eta$, there
exists $b^{\ast}(\eta, l)$ such that $\forall ~ |b|\leq b^{\ast}(\eta, l)$, the estimates of this proposition  hold
with universal constants.
\end{remark}
\begin{remark}The gain is that we now have a Lyapunov function $\Im$ in $H^1$. Remark that in a regime when $\epsilon$ is small compared to $b$ in a certain sense, $\Im \sim \displaystyle{\int}\left|Q_{b}\right|^2 - \displaystyle{\int}\left|Q\right|^2 \sim b^2$ from (\ref{supercritical mass}). Can occur from (\ref{derniere estimation}), this forces $b$ to decay.
\end{remark}
\textbf{Proof:} Multiply (\ref{f1s}) by $\frac{\delta_{1}b}{800}$ and sum with (\ref{control de phiA}). We get 
\begin{align}\label{demo de phiA}
\big(\int\phi_{A}\left|\epsilon\right|^2\big)_{s} + \frac{\delta_{1}b}{800}(f_{1})_{s} &\geq \frac{\delta^{2}_{1}b}{800}\big(\int\left|\nabla \tilde{\epsilon}\right|^2 + \int\left|\tilde \epsilon\right|^2e^{-\left|y\right|}\big) + \frac{b}{800}\int_{A}^{2A}\left|\epsilon\right|^2 \nonumber\\
&+ \frac{c\delta_{1}b}{1000}\Gamma_{b} -C\frac{\lambda^2}{b^2}E(u(t))- \Gamma_{b}^{\frac{l}{2}}\int\left|\nabla\epsilon\right|^2,
\end{align}
We first integrate the left-hand side of (\ref{demo de phiA}) by parts in time:
\begin{align}\label{bf1}
b(f_{1})_{s} &= \bigg(b\tilde{f_{1}}(b) - \int_{0}^{b}\tilde{f_{1}}(v)dv + b\big((\epsilon_{2},(\tilde{\zeta_{re}})_{d}) - (\epsilon_{1},(\tilde{\zeta_{im}})_{d})\big)\bigg)_{s}\nonumber\\
&-b_{s}\big((\epsilon_{2},(\tilde{\zeta_{re}})_{d}) - (\epsilon_{1},(\tilde{\zeta_{im}})_{d})\big),
\end{align}
where $\tilde{f_{1}}$ given by (\ref{f1tilde}). (\ref{bf1}) now yields 
\begin{align}
&\bigg(\int\phi_{A}\left|\epsilon\right|^2  + \frac{\delta_{1}}{800}\big(b\tilde{f_{1}}(b) - \int_{0}^{b}\tilde{f_{1}}(v)dv + b\big((\epsilon_{2},(\tilde{\zeta_{re}})_{d}) - (\epsilon_{1},(\tilde{\zeta_{im}})_{d})\big)\bigg)_{s}\nonumber\\
&\geq \frac{\delta^{2}_{1}b}{800}\big(\int\left|\nabla \tilde{\epsilon}\right|^2 + \int\left|\tilde \epsilon\right|^2e^{-\left|y\right|} + \int_{A}^{2A}\left|\epsilon\right|^2\big) + \frac{c\delta_{1}b}{1000}\Gamma_{b} - C\frac{\lambda^2}{b^2}E(u(t))\nonumber\\
&-\Gamma_{b}^{\frac{l}{2}}\int\left|\nabla \epsilon\right|^2 + \frac{\delta_{1}}{800}b_{s}\big((\epsilon_{2},(\tilde{\zeta_{re}})_{d}) - (\epsilon_{1},(\tilde{\zeta_{im}})_{d})\big).\nonumber
\end{align}
We now inject the expression of the $L^2$-norm:\\
$$\int\left|\epsilon\right|^2 + \int\left|Q_{b}\right|^2 + 2(\epsilon_{1},\Sigma) + 2(\epsilon_{2}, \Theta) = e^{-2at}\int\left|u_{0}\right|^2dx,$$

but $$\int\phi_{A}\left|\epsilon\right|^2 = \int\left|\epsilon\right|^2 - \int(1-\phi_{A})\left|\epsilon\right|^2,$$ we compute:
\begin{align}
\big(\int\phi_{A}\left|\epsilon\right|^2\big)_{s} &= -\bigg(\big(\int\left|Q_{b}\right|^2 - \int\left|Q\right|^2\big) + 2(\epsilon_{1}, \Sigma) + 2(\epsilon_{2}, \Theta) + \int(1-\phi_{A})\left|\epsilon\right|^2\bigg)_{s}\nonumber\\
&+\big(\int\left|u_{0}\right|^2e^{-2at}\big)_{s}\nonumber\\
&=-\bigg(\big(\int\left|Q_{b}\right|^2 - \int\left|Q\right|^2\big) + 2(\epsilon_{1}, \Sigma) + 2(\epsilon_{2}, \Theta) + \int(1-\phi_{A})\left|\epsilon\right|^2\bigg)_{s}\nonumber\\
&-2a\lambda^2e^{-2at}\left\|u_{0}\right\|_{L^2(\mathbb{R}^{d})}^2.\nonumber
\end{align}
Thus, we get
\begin{align}\label{estimation sur Fs}
(-\Im)_{s} &\geq \frac{\delta^{2}_{1}b}{800}\biggl(\int\left|\nabla \tilde{\epsilon}\right|^2 + \int\left|\tilde \epsilon\right|^2e^{-\left|y\right|} + \int_{A}^{2A}\left|\epsilon\right|^2\biggr) + \frac{c\delta_{1}b}{1000}\Gamma_{b} - C\frac{\lambda^2}{b^2}E(u(t))\nonumber\\
&-\Gamma_{b}^{\frac{l}{2}}\int\left|\nabla \epsilon\right|^2 + \frac{\delta_{1}}{800}b_{s}\big((\epsilon_{2},(\tilde{\zeta_{re}})_{d}) - (\epsilon_{1},(\tilde{\zeta_{im}})_{d})\big) -C_{0}\lambda^2.
\end{align}
We now have \\
$$\Gamma_{b}^{\frac{l}{2}}\int\left|\nabla\epsilon\right|^2 \leq \Gamma_{b}^{\frac{l}{2}}\big(\Gamma_{b}^{1-C\eta} + \int\left|\nabla \tilde\epsilon\right|^2\big),$$

from the assumption $l > C\eta$. Next, we estimate from (\ref{control 3}):
\begin{align}
\left|b_{s}\big((\epsilon_{2},(\tilde{\zeta_{re}})_{d}) - (\epsilon_{1},(\tilde{\zeta_{im}})_{d})\big)\right| \leq \Gamma_{b}^{\frac{1}{2} - C\eta}\biggl(\int{\left|\nabla \epsilon\right|^2 + \int\left|\epsilon\right|^2e^{-\left|y\right|} + C\lambda^{2}E(u(t))}\biggr).\nonumber
\end{align}
Injecting these estimates into (\ref{estimation sur Fs}) yields (\ref{derniere estimation}). This concludes the proof of the proposition.
$\hfill\Box$\\\\
Note that now  from (\ref{derniere estimation}) we obtain:
\\
\centerline{$(\Im(s))_{s} \leq -Cb\Gamma_b + C\frac{\lambda^2}{b^2}E(u(t)) \leq -\frac{1}{2}Cb\Gamma_{b} \leq 0$.}

\subsection{Proof of the Bootstrap (Proposition \ref{bootstrap lemma})}
Let $f_{2}$ be defined by:\\
\centerline{$f_{2} = \bigg(\displaystyle{\int\left|Q_{b}\right|^2 - \int Q^2\bigg) - \frac{\delta_{1}}{800}\bigg(b\tilde{f_{1}}(b) - \int_{0}^{b}\tilde{f_{1}}(v)dv}\bigg)$}
\bigskip
it satisfies(using the smallness (\ref{estimation su zeta tilde}) of $\zeta$ in 
$L^2$)
\begin{equation}\label{d0}
\frac{d_{0}}{C} < \frac{df_{2}}{db^2}|_{b^2 =0} < Cd_{0},
\end{equation}
with $d_{0}$ defined by :
\begin{equation}\label{defd0}
0 < \frac{d}{db^2}\big(\int|Q_{b}|^2\big)|_{b^2=0} = d_{0} < +\infty.
\end{equation}
Now as a consequence of the control (\ref{control 1}) and of the coercivity of the linearized energy under the chosen set of orthogonality conditions, we have the bounds:
\begin{equation}\label{bound}
\bigg(\Im(s) - f_{2}(b(s))\bigg) 
\left\{ \begin{array}{lll}
 \geq -\Gamma_{b}^{1 - Cl} + \frac{1}{C}\big(\int\left|\nabla \epsilon\right|^2 + \int\left|\epsilon\right|^2e^{-\left|y\right|}\big)\\
\leq CA^2\big( \int\left|\nabla \epsilon\right|^2 + \int\left|\epsilon\right|^2e^{-\left|y\right|}\big) +\Gamma_{b}^{1 - Cl}.\end{array} \right.
\end{equation} 

Indeed, 
\begin{align}\label{proofdifficile}
\Im(s) - f_{2}(b(s)) & = 2(\epsilon_1,\Sigma) + 2 (\epsilon_2,\Theta) + \int(1-\phi_A)\left|\epsilon\right|^2\nonumber\\
&-\frac{\delta_1b}{800}\bigg(\big(\epsilon_{2},(\tilde{\zeta}_{re})_d\big) - \big(\epsilon_{1}, (\tilde{\zeta}_{im})_d\big)\bigg).
\end{align}

From the estimates on $\tilde{\zeta}$ of Lemma \ref{Linear outgoing radiation}, the choice of $A$, we have:
\begin{align}
 \left|\big(\epsilon_{2},(\tilde{\zeta}_{re})_d\big) - \big(\epsilon_{1}, (\tilde{\zeta}_{im})_d\big)\right|&\leq \Gamma_b^{\frac{1}{2}-C\eta}\big(\int_0^A\left|\epsilon\right|^2\big)^{\frac{1}{2}}\nonumber\\
&\leq A^2\Gamma_b^{1-C\eta} + C(\int\left|\nabla \epsilon\right|^2 + \int\left|\epsilon\right|^2e^{-\left|y\right|})\nonumber\\
&\leq \Gamma_{b}^{1-Cl} + C(\int\left|\nabla \epsilon\right|^2 + \int\left|\epsilon\right|^2e^{-\left|y\right|}).\nonumber
\end{align} 
The other term in (\ref{proofdifficile}) is estimated from the expression of energy:
\begin{align}
&2(\epsilon_{1}, \Sigma + b\Theta_{d} - Re(\Psi)) + 2(\epsilon_{2}, \Theta - b\Sigma_{d} - Im(\Psi))\nonumber\\
&= 2E(Q_{b}) - 2\lambda^2E(u(t))\nonumber\\
&+\int{\left|\nabla\epsilon\right|^2} - \int(\frac{4\Sigma^2}{d\left|Q_{b}\right|^2} + 1)\left|Q_{b}\right|^{\frac{4}{d}}\epsilon_{1}^2\nonumber\\
&-\int(\frac{4\Theta^2}{d\left|Q_{b}\right|^2} + 1)\left|Q_{b}\right|^{\frac{4}{d}}\epsilon_{2}^2\nonumber\\
&-8\int\frac{\Sigma\Theta}{d\left|Q_{b}\right|^2}\left|Q_{b}\right|^{\frac{4}{d}}\epsilon_{1}\epsilon_{2} -\frac{2}{2+\frac{4}{d}}\int F(\epsilon),\nonumber
 \end{align}
which can be rewritten as :
\begin{align}
&2(\epsilon_{1}, \Sigma ) + 2(\epsilon_{2}, \Theta) + \int(1-\phi_A)\left|\epsilon\right|^2 = (L_+\epsilon_1,\epsilon_1) + (L_-\epsilon_2,\epsilon_2) - \int\phi_A\left|\epsilon\right|^2\nonumber\\
& + 2(\epsilon_{1}, Re(\Psi)) + 2(\epsilon_{2}, Im(\Psi)) +2E(Q_{b}) - 2\lambda^2E(u(t))\nonumber\\
&-\int\bigg(\frac{4\Sigma^2}{d\left|Q_{b}\right|^2} + 1)\left|Q_{b}\right|^{\frac{4}{d}}-(\frac{4}{d}+1)Q^{\frac{4}{d}}\bigg)\epsilon_{1}^2\nonumber\\
&-\int\bigg((\frac{4\Theta^2}{d\left|Q_{b}\right|^2} + 1)\left|Q_{b}\right|^{\frac{4}{d}}- Q^{\frac{4}{d}}\bigg)\epsilon_{2}^2\nonumber\\
&-8\int\frac{\Sigma\Theta}{d\left|Q_{b}\right|^2}\left|Q_{b}\right|^{\frac{4}{d}}\epsilon_{1}\epsilon_{2} -\frac{2}{2+\frac{4}{d}}\int F(\epsilon).\nonumber
 \end{align}
We first estimate:
\begin{align}
 \left|(\epsilon_{1},Re(\Psi))\right| &+ \left|(\epsilon_{2},Im(\Psi))\right| + E(Q_{b}) + 2\lambda^2\left|E(u(t))\right|\nonumber\\
&\leq \Gamma_{b}^{1-Cl} + \Gamma_{b}^{l}(\int\left|\nabla \epsilon\right|^2 + \int \left|\epsilon\right|^2e^{-\left|y\right|}).\nonumber
\end{align}
The cubic term $\displaystyle{\int}\left|F(\epsilon)\right|$ and the rest of the quadratic form are controlled by $\delta(\alpha_0) (\displaystyle{\int}\left|\nabla \epsilon\right|^2 + \displaystyle{\int} \left|\epsilon\right|^2e^{-\left|y\right|}) + \Gamma_{b}^{1+l}$,
we thus obtain:
\begin{align}
 &\left|\Im(s) - f_{2}(b(s)) - \bigg((L_{+}\epsilon_1,\epsilon_1) + (L_{-}\epsilon_2,\epsilon_2) - \int(\phi_A)\left|\epsilon\right|^2\bigg)\right|\nonumber\\
&\leq \delta(\alpha_0) (\int \left|\nabla \epsilon\right|^2 + \int \left|\epsilon\right|^2e^{-\left|y\right|}) + \Gamma_{b}^{1-Cl}.\nonumber
\end{align}
The upper bound follows from:
\begin{equation}
 \int(1-\phi_A)\left|\epsilon\right|^2 \leq CA^2\text{log} A(\int\left|\nabla \epsilon\right|^2 + \int \left|\epsilon\right|^2e^{-\left|y\right|}).
\end{equation}
For the lower bound, we use the elliptic estimate on $L = (L_+,L_-)$ (for the proof see Appendix D in \cite{Merle 4}), this ends the proof of (\ref{bound}).\\
We are now in a position to prove the pointwise bound (\ref{control of epsilon 2}):\\
\centerline{$\displaystyle{\int\left|\nabla \epsilon\right|^2 + \int \left|\epsilon\right|^2e^{-\left|y\right|}} \lesssim \Gamma_{b(s)}^{\frac{2}{3}}.$}
Let $s_{2} \in [s_{0},s^{+}]$, if $b_{s}(s_{2}) \leq 0$, then (\ref{control of epsilon 2}) follows directly from (\ref{control de bs}). If $b_{s}(s_{2}) > 0$, let $s_{1} \in [s_{0},s^{+}]$ be � from $s_{2}$ such that $b_{s}(s_{1}) = 0$, then either $s_{1}$ is attained or $s_{1} = s_{0}$. In both cases, we have using (\ref{smallness dans H1}):\\
 \centerline{$\displaystyle{\int\left|\nabla \epsilon(s_{1})\right|^2 + \int \left|\epsilon(s_{1})\right|^2e^{-\left|y\right|}} \leq \Gamma_{b(s_{1})}^{\frac{3}{4}}$}
 and thus 
 \begin{equation}\label{F-f}
 \displaystyle{\Im(s_{1}) - f_{2}(b(s_{1})) \leq \Gamma_{b(s_{1})}^{\frac{17}{24}},}
 \end{equation}
from (\ref{bound}) and for $l>0$ is small enough. Moreover, $b_{s} \geq 0$ on $[s_{1},s_{2}]$ and thus:\\
\begin{equation}\label{bs2-bs1}
\centerline{$
b(s_{2}) \geq b(s_{1})$.}
\end{equation}
We now use the Lyapunov control (\ref{derniere estimation}) to derive:\\
$$\Im(s_{2})\leq \Im(s_{1}),$$
%and thus \\
%$$\Im(s_{2})  \leq \Im(s_{1}) + \lambda^{\alpha}(s_{1}) \leq \Im(s_{1}) + \Gamma_{b(s_{1})}$$
%\bigskip
we then inject (\ref{bound}), (\ref{bs2-bs1}) and (\ref{F-f}) to conclude:
\begin{align}
f_{2}(b(s_{2})) &+ \frac{1}{C}\bigg(\int\left|\nabla \epsilon(s_{2})\right|^2 + \int\left|\epsilon(s_{2})\right|^2e^{-\left|y\right|}\bigg) \leq \Im(s_{2}) + \Gamma_{b(s_{2})}^{1 -Cl}\nonumber\\
&\leq f_2(b(s_1)) + \Gamma_{b(s_{1})}^{1 -Cl} + \Gamma_{b(s_{2})}^{1 -Cl} \leq f_{2}(b(s_{1})) + 2\Gamma_{b(s_{2})}^{\frac{2}{3}}.\nonumber
\end{align} 
The monotonicity (\ref{d0}) of $f_{2}$ in $b$ and (\ref{bs2-bs1}) now imply:\\
\centerline{$\displaystyle{\int\left|\nabla \epsilon(s_{2})\right|^2 + \int \left|\epsilon(s_{2})\right|^2e^{-\left|y\right|}} \lesssim \Gamma_{b(s_{2})}^{\frac{2}{3}}$}
which implies that Equations (\ref{control of epsilon 2}) holds at $s_{2}$. This concludes the proof of (\ref{control of epsilon 2}).\\

Now we are going to prove the upper bound on blowup rate:\\

%Let 
%$$\tilde{b}(s) = b(s) - \int_{s}^{s^{+}}\lambda^{\alpha}(s)ds$$

%then from (\ref{control of lambda}),
%$$\frac{9}{10}b(s) \leq \tilde{b}(s) \leq \frac{10}{9}b(s)$$
From (\ref{control de bs}) we obtain:

\begin{equation}\label{bstilde}
\displaystyle{ b_{s} \geq -\Gamma_{b(s)}^{1 -C\eta}.}
\end{equation}
\\
In particular:
$$ \big(e^{\frac{\pi}{2b(s)}}\big)_{s} \leq e^{\frac{\pi}{2b(s)}}\frac{\pi\Gamma^{1-C\eta}}{2b^2} \leq 1$$
as $\Gamma_{b} \sim e^{-\frac{\pi}{b}}$, and therefore 
$$e^{\frac{\pi}{2b(s)}} \leq e^{\frac{\pi}{b(0)}} + s - s_{0} \leq s,$$
 thus from (\ref{bstilde}) and the value of $s_{0}$(we take $s_{0} = e^{\frac{5\pi}{9b(0)}}$).\\
 Finally, $\forall s \in [s_{0}, s^{+}[$, 
 \begin{equation}\label{bpluslog}
 \displaystyle{b(s) \geq \frac{\pi}{2\log(s)},}
 \end{equation}
 we now rewrite the estimate (\ref{control 3}) using (\ref{control of epsilon 2}) as follows
 \begin{equation}\label{lambdas sur lambda}
 \displaystyle \left|\frac{\lambda_{s}}{\lambda} + b\right| \leq \Gamma_{b}^{\frac{1}{4}}.
 \end{equation}
 Thus:
 
 $$\frac{b}{2} \leq -\frac{\lambda_{s}}{\lambda} \leq 2b.$$
 
 We integrate this in time and get: $\forall s \in [s_{0}, s^{+}]$,
 $$-\log\left(\lambda(s)\right) \geq -\log\lambda\left((s_{0})\right) + \frac{1}{2}\displaystyle{\int_{s_{0}}^{s}b} \geq -\text{log}\lambda(s_{0}) +\frac{\pi}{4}\big(\frac{s}{\text{log}s} - \frac{s_{0}}{\text{log}s_{0}}\big). $$
 Now from (\ref{condition sur b et lambda}):
 $$-\log\left(\lambda(0)\right) \geq e^{\frac{2\pi}{3b(s_{0})}} = s_{0}^{\frac{4}{3}},$$
 and thus 
 \begin{equation}\label{loglambda}
 \displaystyle -\log\left(\lambda(s)\right) \geq -\frac{2}{3}\log\left(\lambda(0)\right) + \frac{\pi}{4}\frac{s}{\log(s)},\quad \text{i.e}\quad \lambda(s) \leq \lambda^{\frac{2}{3}}(0)e^{-\frac{\pi}{4}\frac{s}{logs}}.
 \end{equation}
 This also implies: $\forall s \in [s_{0},s_{2}[$,
 \begin{equation}\label{logetracine}
 \displaystyle{-\text{log}\left(\lambda(s)\right) \geq \frac{\pi}{4}\frac{s}{\text{log(s)}} \geq \sqrt{s}}
 \end{equation}
 and taking the $\log$ of this inequality yields
 \begin{equation}\label{b et loglog}
 \displaystyle{\text{log}\left|\text{log}\left(\lambda(s)\right)\right| \geq \frac{1}{2}\text{log}(s), \quad \text{i.e}\quad b \geq \frac{\pi}{4\text{log}\left|\text{log}\lambda(s)\right|}}
 \end{equation}
 using (\ref{bpluslog}). Therefore (\ref{control of lambda 2}) is proved.\\
 
 Now we are going to prove the monotonicity of $\lambda$:
 We turn to the proof of (\ref{monotonicity 2}) and (\ref{tk}). From (\ref{control of lambda 2}), (\ref{control 3}) and (\ref{control of epsilon 2}), there holds:
 
 $$-\frac{\lambda_{s}}{\lambda} \geq \frac{b}{2}$$

 and thus: $\forall s_{1}$, $s_{2} \in [s_{0},s^{+}]$,
 \begin{equation}\label{logs2surs1}
 \displaystyle{-\text{log}\biggl(\frac{\lambda(s_{2})}{\lambda(s_{1})}\biggr) \geq \frac{1}{2}\int_{s_{1}}^{s_{2}}b  \geq \frac{1}{2}\int_{s_{1}}^{s_{2}}\frac{ds}{\text{log}(s)}.}
 \end{equation}
 This prove (\ref{monotonicity 2}).
 To prove (\ref{tk}), we let $[t_{k},t_{k+1}]$ be a doubling time interval, then from (\ref{logs2surs1}):
 
$$\text{log(2)} = -\text{log}\big(\frac{\lambda(t_{k+1})}{\lambda(t_{k})}\big) \geq \frac{1}{2}\int_{t_{k}}^{t_{k+1}}\frac{dt}{\lambda^2(t)\text{log}(s(t))}$$
 and thus
 $$ 1\geq \frac{C(t_{K+1} - t_{k})}{\lambda^2(t_{k})\text{log}(s(t_{k+1}))} \geq \frac{C(t_{K+1} - t_{k})}{\lambda^2(t_{k})\text{log}\left|\text{log}(\lambda(t_{k+1}))\right|} \geq \frac{C(t_{K+1} - t_{k})}{\lambda^2(t_{k})\text{log}k}$$
 and (\ref{tk}) follows.\\
   The upper bound  on $b$  is a direct consequence (\ref{supercritical mass}). The lower bound $b > 0$ follows
from (\ref{b et loglog}).\\

 Now we prove the blowup in finite time.
 We observe from (\ref{logetracine}) that 
$$T=\displaystyle{\int_{0}^{+\infty}}\lambda^2(s)ds \leq \lambda^{\frac{4}{3}}(0)(C + \displaystyle{\int_{2}^{+\infty}}e^{-\frac{\pi}{4}\frac{s}{logs}}ds) < +\infty.$$
 
 Moreover, from (\ref{control of epsilon 2}),
 
 $$\left\|u(t)\right\|_{H^1(\mathbb{R}^{d})} \sim \frac{1}{\lambda(t)},$$
 and thus the local well-posedness theory in $H^1$ ensures $\lambda(t) \rightarrow 0$ as $t\rightarrow T$.\\
 The convergence of the concentration point is a consequence of (\ref{control 4}), (\ref{logetracine}) and (\ref{control of epsilon 2}) which imply:\\
 \centerline{$\displaystyle{\int_{s_{0}}^{+\infty}}\left|x_{s}\right|ds \leq C \displaystyle{\int_{s_{0}}^{+\infty}}\lambda(s)ds < +\infty.$}
\section{Determination of the blow-up speed}
 In this part, we prove that the   blow-up holds with the log-log speed.\\
 
 Observe from (\ref{bound}), (\ref{control of epsilon 2}) and (\ref{control of lambda 2}) that:\\
 $$\frac{b^2(s)}{C} \leq \Im(s) \leq Cb^2(s).$$
 Together with (\ref{derniere estimation}), this implies:\\
$$\big(\Im\big)_{s} \leq e^{-\frac{C}{\sqrt{\tilde{\Im}}}},$$
 integrating this in time yields:
 $$b(s) \leq C\sqrt{\Im}\leq \frac{C}{\text{log}(s)}$$
 for $s$ large enough. Integrating now (\ref{lambdas sur lambda}) in time, we conclude that 
$$-\text{log}\big(\lambda(s)\big) \leq C \displaystyle{\int_{s_{0}}^{s}}b + C \leq C\frac{s}{\text{log(s)}}$$
 for $s$ large enough, and thus together with (\ref{b et loglog}):
 \begin{equation}\label{blog}
 \displaystyle{\frac{1}{C} \leq b\text{log}\left|\text{log}(\lambda)\right| \leq C.}
 \end{equation}
 Now
\begin{align}\label{lambda2loglog}
-\big(\lambda^2\text{log}\left|\text{log}\lambda\right|\big)_{t} &= -\lambda\lambda_{t}\text{log}\left|\text{log}(\lambda)\right|\bigg(2 + \frac{1}{\left|\text{log}\lambda\right|\text{log}\left|\text{log}\lambda\right|}\bigg)\nonumber\\
&=-\big(\frac{\lambda_{s}}{\lambda} +b\big)\text{log}\left|\text{log}\lambda\right|\bigg(2 + \frac{1}{\left|\text{log}\lambda\right|\text{log}\left|\text{log}\lambda\right|}\bigg)\nonumber\\
&+b\text{log}\left|\text{log}\right|\bigg(2 + \frac{1}{\left|\text{log}\lambda\right|\text{log}\left|\text{log}\lambda\right|}\bigg).
\end{align}

From (\ref{lambdas sur lambda})
\begin{align}
&\int_{t}^{T}\left|\biggl(\frac{\lambda_{s}}{\lambda} + b\biggr)\text{log}|\text{log}\lambda|\right| 
\lesssim \int_{t}^{T}\big(b^2\text{log}\left|\text{log}\lambda\right|)dt.\nonumber\\
\end{align} 
%where we used the almost monotonicity of $\lambda$ in the last steps.
Injecting this into (\ref{lambda2loglog}) integrated from $t$ to $T$ and using (\ref{blog}), we conclude that

$$\frac{T-t}{C} \leq \lambda^2(t)\text{log}\left|\text{log}\lambda(t)\right|\leq C(T-t)$$
\bigskip
from which
\begin{equation}\label{regimeloglog}
\displaystyle{\frac{1}{C}\bigg(\frac{T - t}{\text{log}\left|\text{log}(T-t)\right|}\bigg)^{\frac{1}{2}} \leq \lambda(t) \leq C \bigg(\frac{T - t}{\text{log}\left|\text{log}(T-t)\right|}\bigg)^{\frac{1}{2}}}
\end{equation}
for $t$ close enough to $T$.\\
But 
$$2b(T-t) = \lambda^{2},$$
and thus:\\
$$\frac{1}{C}\frac{1}{\text{log}\left|\text{log}(T-t)\right|} \leq b(t) \leq C \frac{1}{\text{log}\left|\text{log}(T-t)\right|}.$$
From this we obtain:\\
$$\frac{1}{C}\left|\text{log}(T-t)\right|\text{log}\left|\text{log}(T-t)\right| \leq s(t) \leq C \left|\text{log}(T-t)\right|\text{log}\left|\text{log}(T-t)\right|.$$
This prove (\ref{speed}).\qed
\begin{remark}
 Following the ideas in \cite{Merle5}, we can also prove the existence of a 
 $L^2$-profile at blow-up point. More precisely there exists
 $u^{\ast}$ in $L^2(\mathbb{R}^{d})$ ($L^2$-profile) such that:
%a $L^2$-profile at blow up time:\\
%There exist a profile $u^{\ast}$ in $L^2$ such that:

\begin{equation}\label{profile}
u(t,x)-\frac{1}{\lambda^{\frac{d}{2}}(t)}Q_{b(t)}(t,\frac{x-x(t)}{\lambda(t)})e^{i\gamma(t)} \rightarrow u^{\ast}~ \text{in} ~L^2(\mathbb{R}^{d}), ~ t \rightarrow T.
\end{equation}
\end{remark}
 \textbf{Proof of Corollary \ref{colloraireperturbation}}: Let $S(t)$ be the propagator for the linear equation:
\begin{equation}\label{linear}
 i\partial_{t}u + \Delta u  = 0, \quad (t,x) \in [0,\infty[\times\mathbb{R}^{d}.\nonumber
\end{equation}
The Cauchy problem for (\ref{NLSa}) with $u(0) = u_{0} \in H^{1}(\mathbb{R}^{d})$ is equivalent to the integral equation:
\begin{equation}
 u(t) = S(t)u_{0} + i\int_{0}^{t}S(t-s)(|u(s)|^{\frac{4}{d}}u(s) + iau)ds. \nonumber
\end{equation}
We know from Lemma \ref{us1}: there exist $T(\left\|u_0\right\|_{H^1(\mathbb{R}^{d})}) > 0$ such that:\\
\centerline{ $ \forall 0 \leq a \leq 1$, $\left\|u\right\|_{L^{\infty}([0,T]; H^1)} \leq 2\left\|u_0\right\|_{H^1(\mathbb{R}^{d})}.$}
Let $u$ a solution for (\ref{NLSa}) and $v$ solution for (\ref{NLS}) we have:\\
\begin{align}
u-v = S(t)(u_{0} - v_{0})  &+i \int_{0}^{t}S(t-s)(|u(s)|^{\frac{4}{d}}u(s)-|v(s)|^{\frac{4}{d}}v(s))ds\nonumber\\ &+i a\int_{0}^{t}S(t-t^\prime)u(t^\prime)dt^\prime.\nonumber
\end{align}
By Strichartz we obtain (see the proof of Lemma \ref{us1}):
\begin{align}
\left\|u-v\right\|_{L^{\infty}([0,T]; H^1)} &\leq
\left\|u_{0}-v_{0}\right\|_{H^1(\mathbb{R}^{d})}\nonumber\\ 
&+ C T^{\gamma}\big(
\left\|u\right\|^{p}_{L^{\infty}([0,T]; H^1)} 
+\left\|v\right\|^{p}_{L^{\infty}([0,T]; H^1)}\big)\left\|u-v\right\|_{L^{\infty}([0,T]; H^1)}\nonumber\\
&+CaT \left\|u\right\|_{L^{\infty}([0,T]; H^1)}.\nonumber
\end{align}
Thus for $T_1 = \text{Min}\big(T(\left\|u_0\right\|_{H^1(\mathbb{R}^{d})}),T(\left\|v_0\right\|_{H^1(\mathbb{R}^{d})})\big)$ we obtain $\forall~ 0\leq t \leq T_{1}$:
\begin{align}
\left\|u-v\right\|_{L^{\infty}([0,T]; H^1)} &\leq
\left\|u_{0}-v_{0}\right\|_{H^1(\mathbb{R}^{d})}\nonumber\\ 
&+ C T^{\gamma}\big(
\left\|u_0\right\|^{p}_{ H^1} 
+\left\|v_0\right\|^{p}_{ H^1}\big)\left\|u-v\right\|_{L^{\infty}([0,T]; H^1)}\nonumber\\
&+CaT \left\|u\right\|_{L^{\infty}([0,T]; H^1)}.\nonumber
\end{align} 
Now for $T_2 = \frac{1}{2}\text{Min}\big(\text{Max}^{-\frac{1}{\gamma}}(\left\|u_{0}\right\|^{p}_{ H^1}, \left\|v_{0}\right\|^{p}_{ H^1}), T_1\big)$:

\begin{align}
 \left\|u-v\right\|_{L^{\infty}([0,t]; H^1)} &\leq
\left\|u_{0}-v_{0}\right\|_{H^1(\mathbb{R}^{d})} + \frac{1}{2} \left\|u-v\right\|_{L^{\infty}([0,t]; H^1)}\nonumber\\
& + a,\nonumber
\end{align}
thus\\ 
\centerline{$\left\|u-v\right\|_{L^{\infty}([0,t]; H^1)} \lesssim
\left\|u_{0}-v_{0}\right\|_{H^1(\mathbb{R}^{d})} + a.~\forall ~ 0 < t < T_2$.}
Thus the map $(a, \phi) \rightarrow u(\cdot, a, \phi)$ is continuous in $(0,u_0)$ from $\mathbb{R} \times H^1(\mathbb{R}^d)$ to $C([0,T_2], H^1(\mathbb{R}^{d}))$. Since $T_2$  only depends on $\left\|u_0\right\|_{H^1(\mathbb{R}^{d})}$, this continuity extends to any interval $[0,T]$ in the maximal interval of existence of $u$. \\
%The energy $E$ and the kinetic momentum $P$ are continuous from $H^{1}(\mathbb{R}^{d})$ to $\mathbb{R}$ thus $v_{0}$ verifie the conditions C.I, and by theorem \ref{theoremessentiel} we obtain the blow up of $v$ with the initial data $v_{0}$ close to $u_{0}$ in $H^{1}(\mathbb{R}^{d})$ and $a$ small in $\mathbb{R}$. 
We  know after a time $t_0$ closed to blow-up time of $u$ with the initial data $u_0$, that $ u(t_0)$ verifies C.I, and by continuity  $v(t_0)$ verifies also C.I (the conditions C.I are stable by a small perturbations in $H^1$), then  we obtain from Theorem \ref{theorem 3} the blow up of $v$ with the initial data $v(t_0)$ for (\ref{NLSa}).Therefore the solution of (\ref{NLSa}) emanating from $v_0$ blows up in finite time in the log-log regime.
\hfill$\Box$

%\bibliographystyle{plain}
%\bibliography{biblio1}

\begin{thebibliography}{10}

\bibitem{Ber}
H.~Berestycki and P.-L. Lions.
\newblock{\it Nonlinear scalar field equations. {II}. {E}xistence of infinitely
  many solutions}.
\newblock { Arch. Rational Mech. Anal.}, 82(1983):347--375.

\bibitem{Cazenave}
T.~Cazenave.
\newblock {\it Semilinear {S}chr\"odinger equations}, volume~10 of {\em Courant
  Lecture Notes in Mathematics}.
\newblock New York University Courant Institute of Mathematical Sciences, New
  York, 2003.

\bibitem{Colliander}
J.~Colliander and P.~Raphael.
\newblock {\it Rough blowup solutions to the {$L^2$} critical {NLS}}.
\newblock { Math. Ann.}, 345(2009):307--366.

\bibitem{Fibich}
G.~Fibich and F.~Merle.
\newblock {\it Self-focusing on bounded domains}.
\newblock { Phys. D}, 155(2001):132--158.

\bibitem{Hmidi}
T.~Hmidi and S.~Keraani.
\newblock{\it Blowup theory for the critical nonlinear {S}chr\"odinger equations
  revisited}.
\newblock {Int. Math. Res. Not.}, 46(2005):2815--2828.

\bibitem{Kato}
T.~Kato.
\newblock {\it On nonlinear {S}chr\"odinger equations}
\newblock {Ann. Inst. H. Poincar\'e Phys. Th\'eor.}, 46(1987):113--129.

\bibitem{Kwong1}
M.K Kwong.
\newblock {\it Uniqueness of positive solutions of {$\Delta u-u+u^p=0$} in {${\bf
  R}^n$}}.
\newblock { Arch. Rational Mech. Anal.}, 105(1989):243--266.

\bibitem{Lions1}
P.-L. Lions.
\newblock {\it The concentration-compactness principle in the calculus of
  variations. {T}he locally compact case. {II}}.
\newblock { Ann. Inst. H. Poincar\'e Anal. Non Lin\'eaire}, 1(1984):223--283.

\bibitem{MerleRaphael1}
F.~Merle and P.~Raphael.
\newblock {\it Blow up dynamic and upper bound on the blow up rate for critical
  nonlinear {S}chr\"odinger equation}.
\newblock In { Journ\'ees ``\'{E}quations aux {D}\'eriv\'ees {P}artielles''
  ({F}orges-les-{E}aux, 2002)}, pages Exp. No. XII, 5. Univ. Nantes, Nantes,
  2002.

\bibitem{Merle2}
F.~Merle and P.~Raphael.
\newblock {\it Sharp upper bound on the blow-up rate for the critical nonlinear
  {S}chr\"odinger equation}.
\newblock { Geom. Funct. Anal.}, 13(2003):591--642.

\bibitem{Merle3}
F.~Merle and P.~Raphael.
\newblock {\it On universality of blow-up profile for {$L^2$} critical nonlinear
  {S}chr\"odinger equation}.
\newblock { Invent. Math.}, 156(2004):565--672.

\bibitem{Merle5}
F.~Merle and P.~Raphael.
\newblock {\it Profiles and quantization of the blow up mass for critical nonlinear
  {S}chr\"odinger equation}.
\newblock {Comm. Math. Phys.}, 253(2005):675--704.

\bibitem{Merle4}
F.~Merle and P.~Raphael.
\newblock {\it On a sharp lower bound on the blow-up rate for the {$L^2$} critical
  nonlinear {S}chr\"odinger equation}.
\newblock { J. Amer. Math. Soc.}, 19(2006):37--90 (electronic).

\bibitem{ohta}
M.~Ohta and G.~Todorova.
\newblock {\it Remarks on global existence and blowup for damped nonlinear
  {S}chr\"odinger equations}.
\newblock { Discrete Contin. Dyn. Syst.}, 23(2009):1313--1325.

\bibitem{Planchon1}
F.~Planchon and P.~Rapha{\"e}l.
\newblock {\it Existence and stability of the log-log blow-up dynamics for the
  {$L^2$}-critical nonlinear {S}chr\"odinger equation in a domain}.
\newblock {Ann. Henri Poincar\'e}, 8(2007):1177--1219.

\bibitem{Merle6}
P.~Raphael.
\newblock {\it Stability of the log-log bound for blow up solutions to the critical
  non linear {S}chr\"odinger equation}.
\newblock { Math. Ann.}, 331(2005):577--609.

\bibitem{Tsutsumi}
M.~Tsutsumi.
\newblock {\it Nonexistence of global solutions to the {C}auchy problem for the
  damped nonlinear {S}chr\"odinger equations}.
\newblock { SIAM J. Math. Anal.}, 15(1984):357--366.

\bibitem{weinst}
M.I. Weinstein.
\newblock {\it Nonlinear {S}chr\"odinger equations and sharp interpolation
  estimates}.
\newblock { Comm. Math. Phys.}, 87(1982/83):567--576.

\end{thebibliography}

%%%%%%%%%%%%%%%%%%%%%%%%%%%%%%%%%%%%%%%%%%%%%%

\end{document}